%% file: Psdodiff_models_US_waves_v1.tex
\documentclass[hidelinks,onefignum,onetabnum]{siamart220329}

\input{ex_shared}

\ifpdf
\hypersetup{
  pdftitle={Pseudodifferential Models for Ultrasound Waves with Fractional Attenuation},
  pdfauthor={Sebastian Acosta, Jesse Chan, Raven Johnson and Benjamin Palacios}
}
\fi

\begin{document}

\maketitle

\begin{abstract}
To strike a balance between modeling accuracy and computational efficiency for simulations of ultrasound waves in soft tissues, we derive a pseudodifferential factorization of the wave operator with fractional attenuation. This factorization allows us to approximately solve the Helmholtz equation via one-way (transmission) or two-way (transmission and reflection) sweeping schemes tailored to high-frequency wave fields. We provide explicitly the three highest order terms of the pseudodifferential expansion to incorporate the well-known square-root first order symbol for wave propagation, the zeroth order symbol for amplitude modulation due to changes in wave speed and damping, and the next symbol to model fractional attenuation. We also propose wide-angle Pad\'e approximations for the pseudodifferential operators corresponding to these three highest order symbols. Our analysis provides insights regarding the role played by the frequency and the Pad\'e approximations in the estimation of error bounds. We also provide a proof-of-concept numerical implementation of the proposed method and test the error estimates numerically. 
\end{abstract}

\begin{keywords}
Wave propagation, acoustics, high frequency, Helmholtz equation, pseudodifferential calculus, Pade approximations
\end{keywords}

\begin{MSCcodes}
35S05, 35S15, 35L05, 41A21, 41A28
\end{MSCcodes}

\section{Introduction}

The present study is motivated by the application of ultrasound waves to medical therapeutics and diagnosis where mathematical modeling and computational simulations are playing an increasingly important role \cite{Ammari-2008,AmmariBook2009,Beard2011,Soldati2014,Verweij2014,AmmariBook2012,Poudel2019,Li2023,Li2023a}. To unleash the full potential of computer simulations of ultrasound fields, the underlying mathematical model of wave propagation must strike the right balance between accuracy and speed of computation. Most methods for ultrasound wave simulation lie at one of two extremes. At one extreme, closed-form models provide tremendous computational speed but sacrifice too much accuracy.
These over-simplified methods are only valid under stringent assumptions and render a multitude of inaccuracies for realistic biological media 
\cite{Webb2003,Verweij2014,Opielinski2000,Guasch2020,Shortell2017,Wiskin2012,Li2010}. At the other extreme, full-waveform models based on partial differential equations are very accurate but suffer serious computational scalability limitations, especially at ultrasonic MHz frequencies. In addition to short wavelengths, ultrasound waves in biological media are characterized by their fractional or power-law attenuation profiles which further complicates the validity of closed-form methods and the scalability of full-waveform simulations \cite{Nasholm2011,Acosta-Montalto-2016,Acosta-Palacios-2018,Kowar2011a,Chen2016,Holm2011,Lou2017,Huang2012a}. This challenge (recognized by recently several research groups \cite{Verweij2014,Gu2018,Leung2019,Gu2019,Haqshenas2021} including the benchmark article recently published by Aubry et al. \cite{Aubry2022}) exposes an unresolved need to simulate ultrasound fields efficiently enough for some clinical environments. 

In most ultrasound medical applications, the transducer is designed to emit waves moving along a dominant direction of propagation. Such is the case for plane waves and beams. Some processing methods (such as ray and wavefront tracing, one-way, parabolic, paraxial, kinetic, and Eikonal models) rely on this geometry of propagation to reduce computational costs \cite{Halpern1988,Bamberger1988b,Engquist2003,Stolk2004,Varslot2005,Engquist2009book,Tanushev2009,Jing2011,OptRoot2010,Angus2014,Gu2018}. In the present work, we pursue a similar line of work. We propose a sweeping algorithm based on a pseudodifferential factorization of the wave equation to improve the accuracy of closed-form methods while keeping the execution times much more competitive than full-waveform simulations. The pseudodifferential calculus allows us
to take advantage of the geometry of acoustic energy flow, incorporate physical interactions such as transmission, reflection and attenuation while efficiently handling the highly oscillatory nature of ultrasound waves. We extend related approaches by Zhang, Stolk, Op't Root, Angus, Vasilyev and others \cite{Zhang2003,Stolk2004,OptRoot2010,Angus2014,Vasilyev2018} by considering the pseudodifferential calculus to all orders, implementing a two-way model to account for reflections, and deriving the pseudodifferential symbol for fractional attenuation. This latter feature is important for applications in biomedicine where fractional or power-law attenuation models have been shown to represent biological media accurately \cite{Treeby2010,Acosta-Palacios-2018,Holm2011,AmmariBook2012,Chen2016}. Our work can also be understood as an alternative to the Bremmer series which has been used for the theoretical analysis and computational treatment of multi-dimensional inverse scattering problems \cite{Shehadeh2017,Malcolm2005}.

This paper is structured as follows. In Section \ref{Section.ReviewPDO} we briefly review the mathematical preliminaries of pseudodifferential calculus. In Section \ref{Section.PSD}, we develop the pseudodifferential factorization of the wave equation with a fraction attenuation term, carry out the calculations for the classical decomposition of the pseudodifferential operators in inverse powers of the frequency $\omega$, and propose high order (or so-called wide-angle) Pad\'e approximations for the pseudodifferential operators corresponding to the highest degree symbols. This development includes the well-known square-root symbol to model one-way wave propagation, the zeroth degree symbol for amplitude modulation due to changes in wave speed and damping, and the next symbol to model fractional attenuation. In Section \ref{Section.Sweeping}, we incorporate the approximated symbols into the proposed two-way sweeping methods. Section \ref{Section.ErrorAnalysis} contains basic error analysis at the continuous level and provides insights concerning the role played by the frequency $\omega$ and the Pad\'e approximations in the estimation of error bounds. 
In Section \ref{Section.Numerical} we provide a proof-of-concept numerical implementation of the proposed method and test the error estimates numerically. Finally, in Section \ref{Section.Conclusion} we offer some concluding remarks, discuss limitations of the proposed pseudodifferential method and areas of potential improvement.

\section{Preliminaries}
\label{Section.ReviewPDO}

We introduce some notation in this section which is self-contained and is not to be confused with the notation in the remainder of the paper. We briefly introduce the definition of pseudodifferential operators and the basic properties of pseudodifferential calculus. We use \cite{Taylor1991Book,Grigis1994Book} as our references. The first ingredient in this formulation is the $n$-dimensional Fourier transform $\mathcal{F}$ and its inverse $\mathcal{F}^{-1}$, which for an admissible function $u: \mathbb{R}^n \to \mathbb{C}$, are respectively given by
\begin{align} \label{Eqn.FT}
\mathcal{F} u(\xi) = \hat{u}(\xi) = (2 \pi)^{-d/2} \int e^{-i x \cdot \xi}  u(x)  dx, 
\end{align}
and
\begin{align} \label{Eqn.IFT}
\mathcal{F}^{-1} \hat{u}(x) = (2 \pi)^{-d/2} \int e^{i x \cdot \xi} \hat{u}(\xi)  d\xi.
\end{align}
The Fourier transform satisfies the following properties in relation to differentiation,
\begin{align} \label{Eqn.DTDiff}
D^{\alpha} u (x) = (2 \pi)^{-d/2} \int \xi^\alpha \hat{u}(\xi) e^{i x \cdot \xi} d\xi
\end{align}
for a multi-index $\alpha=(\alpha_1, ..., \alpha_n) \in \mathbb{N}^n$ where $D^{\alpha} = D_{1}^{\alpha_1} \dots D_{n}^{\alpha_n}$ and $D_{j} = - i \partial_{x_j}$. A differential operator of order $m$ with variable coefficients $a_{\alpha}=a_{\alpha}(x)$ can be expressed as
\begin{align} \label{Eqn.PDO}
A(x,D) = \sum_{|\alpha| \leq m} a_{\alpha}(x) D^{\alpha}
\end{align}
where $|\alpha| = \alpha_1 + ... + \alpha_n$. Then we have
\begin{align} \label{Eqn.PDO2}
A(x,D) u(x) = (2 \pi)^{-d/2} \int e^{i x \cdot \xi} a(x, \xi) \hat{u}(\xi) d \xi 
\end{align}
where
\begin{align} \label{Eqn.PDO3}
 a(x, \xi) = \sum_{|\alpha| \leq m} a_{\alpha}(x)\xi^\alpha
\end{align}
where $\xi^\alpha =  \xi_{1}^{\alpha_1} \, \xi_{2}^{\alpha_2} \, ... \, \xi_{n}^{\alpha_n}$. We call $a(x,\xi)$ the (full) symbol of the differential operator. The Fourier integral representation \eqref{Eqn.PDO2} of differential operators can be used to generalize them to a larger class known as pseudodifferential operators. Note that for differential operators, the functions $a(x,\xi)$ are polynomial with respect to $\xi$. For pseudodifferential operators, we let these functions belong to larger sets, known as Hormander's symbol classes defined as follows. Take $m \in \mathbb{R}$ and define the symbol class $\SS^m = \SS^m(\mathbb{R}^n \times \mathbb{R}^n)$ to consist of $C^{\infty}$ functions $a(x,\xi)$ satisfying
\begin{align} \label{Eqn.Hclass}
| D_{x}^{\beta} D_{\xi}^{\alpha} a(x,\xi)| \leq C_{\alpha \beta} (1 + |\xi|^2)^{(m - |\alpha|)/2}
\end{align}
for all multi-indices $\alpha$ and $\beta$, and some constants $C_{\alpha \beta}$. Once a symbols $a(x,\xi)$ belongs to a class $\SS^m$, the associated operator $\Op{a}$, defined by
\begin{align} \label{Eqn.GeneralPSDO}
\Op{a} u(x) = (2\pi)^{-d/2} \int a(x,\xi) e^{ix\cdot \xi} \hat{u}(\xi)d\xi,
\end{align}
is said to be a pseudodifferential operator that belongs to $\Op{\SS^m}$. In other words, pseudodifferential operators are defined by their symbols. The smallest possible $m$ that allows a symbol $a(x,\xi)$ to satisfy \eqref{Eqn.Hclass} is known as the order of the pseudodifferential operator.
Note that by \eqref{Eqn.Hclass}, when we take a derivative of a symbol $a(x,\xi)$ with respect to $\xi$, we simply obtain another pseudodifferential operator with lower order.

Given a pseudodifferential operator $A\in \Op{\SS^m}$, its symbol will be denoted $\Sym{A}$. An important subset of pseudodifferential operators are those that satisfy the following conditions. Let $a(x,\xi)=\Sym{A} \in \SS^m$. If there are smooth functions $a_{m-j}(x,\xi)$ positively homogeneous in $\xi$ of degree $m-j$, i.e., $a_{m-j}(x,r \xi) = r^{m-j} a_{m-j}(x,\xi)$ for any $r>0$ and $\xi \neq 0$, and if
\begin{align} \label{Eqn.ClassicalSymbol}
\left( a(x,\xi) - \sum_{j=0}^{J} a_{m-j}(x,\xi)  \right) \in \SS^{m-J-1}
\end{align}
for all $J=0, 1, 2, ...$, we say that $A$ is a ``classical" pseudodifferential operator. In such cases when $A$ is classical, we may write
\begin{align} \label{Eqn.ClassicalSymbol2}
a \sim \sum_{j=0}^{\infty} a_{m-j}  \quad \text{or} \quad \Op{a} \sim \sum_{j=0}^{\infty} \Op{a_{m-j}},
\end{align}
and say that $\sum_{j=0}^{\infty} a_{m-j}$ is an asymptotic sum of $a$.

Given two symbols $a(x,\xi)\in \SS^{m_1}$ and  $b(x,\xi)\in \SS^{m_2}$, the composition of their respective operators has a symbol in $\SS^{m_1+m_2}$ satisfying
\begin{align} \label{Eqn.Symbol_of_Composition}
\Sym{\Op{a}\Op{b}} \sim  \sum_{\alpha\geq 0}\frac{i^{|\alpha|}}{\alpha!}D^\alpha_\xi a \, D^\alpha_x b.
\end{align} 

In particular, the operator associated with the product of two symbols $a$ and $b$ satisfies,
\begin{align} \label{Eqn.Operator_of_Product}
\Op{ab} \sim \Op{a} \Op{b} - \sum_{\alpha\geq 1}\frac{i^{|\alpha|}}{\alpha!} \Op{D^\alpha_\xi a \,  D^\alpha_x b }.
\end{align}
Hence if $a = a(x)$ is independent of $\xi$ or $b=b(\xi)$ is independent of $x$, then $\Op{ab} \sim \Op{a} \Op{b}$. 

Another interesting case is the computation of the operator of a symbol of the form $1/b$ for $b \in \SS^{m}$ being an elliptic symbol. This is a special case of the above equations. So we have
\begin{align} \label{Eqn.Operator_of_Inverse}
\Op{b} \Op{1/b} \sim \Op{1} + \sum_{\alpha\geq 1}\frac{i^{|\alpha|}}{\alpha!} \Op{D^\alpha_\xi b \, D^\alpha_x b^{-1} }.
\end{align}
Hence, provided that $\Op{b}$ is invertible, then 
\begin{align} \label{Eqn.Operator_of_Inverse2}
\Op{1/b} \sim \Op{b}^{-1} \left( I + \sum_{\alpha\geq 1}\frac{i^{|\alpha|}}{\alpha!} \Op{D^\alpha_\xi b \,  D^\alpha_x b^{-1} } \right).
\end{align}

Finally, given a pseudodifferential symbol $a \in \SS^{m}$, the Sobolev norms satisfy
\begin{align} \label{Eqn.SobNormPsiOp}
\| \Op{a} u \|^2_{H^{s}(\mathbb{R}^n)} &\leq \int 
 (1+|\xi|^2)^s \sup_{x}{|a(x,\xi)|^2} |\hat{u}(\xi)|^2 d \xi \nonumber \\ 
&\leq C_{00} \int 
 (1+|\xi|^2)^{s+m} |\hat{u}(\xi)|^2 d \xi \nonumber \\
&= C_{00} \| u \|^2_{H^{s+m}(\mathbb{R}^n)} 
\end{align}
provided that $u \in H^{s+m}(\mathbb{R}^n)$ and where $C_{00} > 0$ is the constant appearing in \eqref{Eqn.Hclass}. In other words, $\Op{a}$ maps $H^{s+m}(\mathbb{R}^n)$ continuously into $H^{s}(\mathbb{R}^n)$ for any $s \in \mathbb{R}$. Moreover, the norm of $\Op{a}$ is proportional to the bound on its symbol.

\section{Pseudodifferential factorization of the wave equation}
\label{Section.PSD}

Our approach has evolved from our experience with the formulation and implementation of absorbing boundary conditions for waves \cite{Guo2020,Acosta2021d,Khajah2019,Acosta2017c,Villamizar2017b,Acosta2015f}. 
The starting point is Nirenberg's factorization theorem for hyperbolic differential operators \cite{Nirenberg1973,Antoine1999,Antoine2001b,Acosta2017c}. We consider the wave operator,
\begin{align} \label{Eqn.LaplacianLocal}
\LL u &= \Delta u - c^{-2} \partial_{t}^2 u -a \partial_{t}u- a_\alpha \partial_{t}^{\alpha}u \nonumber \\
 &=  \partial_{x}^2 u  +  \Delta_{\perp} u - c^{-2} \partial_{t}^2 u - a \partial_{t}u - a_{\alpha} \partial_{t}^{\alpha}.
\end{align}
Here, $c$ is the wave speed, $a$ is the damping coefficient, $a_{\alpha}$ is the attenuation coefficient for the fractional term, and $0 < \alpha < 1$, is the fractional exponent of the attenuation.  The solution $u=u(t,x,x_{\perp})$ depends on time $t \in \mathbb{R}$, the coordinate $x \in \mathbb{R}$ along the dominant direction of acoustic propagation, and the coordinates $x_{\perp} \in \mathbb{R}^{d-1}$ perpendicular to the $x$-axis. The spatial dimension of the problem is $d \geq 2$. The Laplacian $\Delta = \partial_{x}^2 + \Delta_{\perp}$ is decomposed accordingly so that $\Delta_{\perp}$ is the tangential Laplacian defined on the hyperplane perpendicular to the $x$-axis. Following Nirenberg's approach, we constructively show the wave operator $\LL$ can be decomposed into forward and backward components based on pseudodifferential operators $\Lambda^{\pm} \in \Op{\SS^1(\mathbb{R}^{d} \times \mathbb{R}^{d})}$, parameterized by $x \in \mathbb{R}$, with respective symbols $\lambda^{\pm}$, such that  
\begin{align} \label{Eqn.Decomp01}
\LL  &= ( \partial_{x} - \Lambda^{-})( \partial_{x} - \Lambda^{+}) =  \partial_{x}^2  -  \left( \Lambda^{+} + \Lambda^{-} \right) \partial_{x}   + \left( \Lambda^{-}\Lambda^{+} - \Op{\partial_{x} \lambda^{+}} \right)
\end{align}
modulo $\Op{\SS^{-\infty}}$, where $ \Op{\partial_{x} \lambda^{+}}$ stands for the pseudodifferential operator with full symbol $\partial_x \lambda_+$. The two families of pseudodifferential symbols $\lambda^{\pm}$ are parameterized by $x \in \mathbb{R}$ whose Fourier dual is $\sigma \in \mathbb{R}$, such that for each fixed $x$, the symbol $\lambda^{\pm} = \lambda^{\pm}(t,x_{\perp},\omega,\sigma_{\perp})$. Here and in what follows, $\omega \in \mathbb{R}$ is the Fourier dual of time $t \in \mathbb{R}$, so that $-\omega^2$ is the symbol of $\partial_{t}^2$. Similarly, $\sigma_{\perp} \in \mathbb{R}^{d-1}$ is the Fourier dual of $x_{\perp} \in \mathbb{R}^{d-1}$, so that $- |\sigma_{\perp}|^2$ is the symbol of $\Delta_{\perp}$. 

We refer to $\Lambda^{+}$ and $\Lambda^{-}$ as the forward and backward Dirichlet--to--Neumann (DtN) operators. Then, by matching terms with same number of $x$-derivatives in (\ref{Eqn.LaplacianLocal}) and (\ref{Eqn.Decomp01}) we obtain,
\begin{align} 
\Lambda^{+} + \Lambda^{-} &= 0,  \label{Eqn.OpMain2} \\
\Lambda^{-}\Lambda^{+} -  \Op{\partial_{x} \lambda^{+}} &= \Delta_{\perp} - c^{-2} \partial_{t}^2  - a \partial_{t} - a_{\alpha} \partial_{t}^{\alpha}.   \label{Eqn.OpMain3}
\end{align}
To process the above equations, we need to obtain the symbol for the product of two pseudodifferential operators as reviewed in Section \ref{Section.ReviewPDO},
\begin{align} 
\Sym{ \Lambda^- \Lambda^+ } &\sim \sum_{m=0}^{\infty} \frac{(-i)^m}{m !} \partial_{\sigma_{\perp}}^m \lambda^{-} \, \partial_{x_{\perp}}^m \lambda^{+}.  \label{Eqn.Prod} 
\end{align}
Note that terms of the form $\partial_{\omega}^m \lambda^{-} \partial_{t}^m \lambda^{+}$ do not appear in the above expression because $\lambda^{\pm}$ are independent of time $t$ since the wave operator \eqref{Eqn.LaplacianLocal} has time-independent coefficients.

The symbols $\lambda^{\pm}$ admit the following classical pseudodifferential expansion 
\begin{align} 
\lambda^{\pm} \sim \sum_{n=-1}^{+ \infty} \lambda^{\pm}_{-n} + \sum_{n=1}^{+ \infty} \lambda^{\pm}_{\beta_{n}} \label{Eqn.SymbolExpansion01}
\end{align}
where $\lambda^{\pm}_{-n} \in \SS^{-n}$ are homogeneous functions of degree $-n$ in $(\omega, \sigma_{\perp})$. Similarly, $\lambda^{\pm}_{\beta_{n}} \in \SS^{\beta_{n}}$ are homogeneous functions of degree $\beta_{n}$ in $(\omega, \sigma_{\perp})$. It is also required that $\{ \beta_{n} \}$ is a strictly decreasing sequence, $\beta_{n} \to -\infty$ as $n \to \infty$, and that $\beta_{n} \notin \mathbb{Z}$ \cite{Antoine1999,Taylor1991Book,Taylor1996-ChapterPSO}.

Plugging \eqref{Eqn.SymbolExpansion01} into the symbolic version of \eqref{Eqn.OpMain2} and collecting symbols of the same degree, we obtain
\begin{align} 
\lambda_{-n}^{+} + \lambda_{-n}^{-} &= 0, \qquad \text{for all $n=-1,0,1,2, ...$},  \label{Eqn.010} \\
\lambda_{\beta_{n}}^{+} + \lambda_{\beta_{n}}^{-} &= 0, \qquad \text{for all $n=1,2, ...$}.  \label{Eqn.011}
\end{align}
Similarly, plugging \eqref{Eqn.SymbolExpansion01} into the symbolic version of \eqref{Eqn.OpMain3}, using \eqref{Eqn.Prod} and collecting terms of the same degree, we get
\begin{align} 
\lambda_{1}^{-} \lambda_{1}^{+} &= \omega^2 c^{-2} - |\sigma_{\perp}|^2 \label{Eqn.020} \\
\lambda_{1}^{-} \lambda_{0}^{+} + \lambda_{0}^{-} \lambda_{1}^{+}  &= \partial_{x} \lambda_{1}^{+} - a i \omega + i  \partial_{\sigma_{\perp}} \lambda_{1}^{-} \, \partial_{x_{\perp}} \lambda_{1}^{+}   \label{Eqn.021} \\
\lambda_{1}^{-} \lambda_{\beta_{1}}^{+} + \lambda_{\beta_{1}}^{-} \lambda_{1}^{+}  &= - a_{\alpha} (i\omega)^{\alpha} \label{Eqn.022}
\end{align}
where $\beta_{1} = \alpha - 1$. Hence, combining \eqref{Eqn.010}-\eqref{Eqn.011} and \eqref{Eqn.020}-\eqref{Eqn.022} we obtain the highest degree symbols
\begin{align} 
\lambda_{1}^{\pm} &= \mp i \sqrt{ \omega^2 c^{-2} - |\sigma_{\perp}|^2 } \label{Eqn.Sym1} \\
\lambda_{0}^{\pm} &= \mp \frac{\partial_{x} \lambda_{1}^+}{ 2 \lambda_{1}^{+} } \pm \frac{a i \omega}{2 \lambda_{1}^{+}} \mp  i \frac{  \partial_{\sigma_{\perp}} \lambda_{1}^{+}  \, \partial_{x_{\perp}} \lambda_{1}^{+} }{2 \lambda_{1}^{+} } \label{Eqn.Sym0} \\
\lambda_{\beta_{1}}^{\pm} &= \pm \frac{a_{\alpha} (i \omega )^{\alpha}}{2 \lambda_{1}^{+}}, \, \text{where $\beta_{1} = \alpha - 1$}. \label{Eqn.SymB0}
\end{align}

For the recursive terms of negative integer order, we collect terms of order $-n$ from \eqref{Eqn.OpMain3} to obtain
\begin{align} 
&\lambda_{1}^{-} \lambda_{-n-1}^{+} + \lambda_{-n-1}^{-} \lambda_{1}^{+} + \sum_{j=0}^{n} \lambda_{-j}^{-} \lambda_{j-n}^{+} 
 + \sum_{m=1}^{n+2} \sum_{j=-1}^{n-m+1}  \frac{(-i)^m}{m!} \partial_{\sigma_{\perp}}^{m}  \lambda_{-j}^{-} \partial_{x_{\perp}}^{m} \lambda_{j-n+m}^{+} \nonumber \\
& + \sum_{\substack{ j,l \\ \beta_{j} + \beta_{l} = - n } } \lambda_{\beta_j}^{-} \lambda_{\beta_l}^{+} + \sum_{m=1}^{n+2} \sum_{\substack{ j,l \\ \beta_{j} + \beta_{l} = - n }}  \frac{(-i)^m}{m!} \partial_{\sigma_{\perp}}^{m}  \lambda_{\beta_{j}}^{-} \partial_{x_{\perp}}^{m} \lambda_{\beta_{l-m}}^{+}  = \partial_{x} \lambda_{-n}^{+} \nonumber
\end{align}
and using \eqref{Eqn.010} we arrive at
\begin{align} 
\lambda_{-n-1}^{\pm} &= \mp  \frac{1}{2 \lambda_{1}^{+} }  \Bigg( \partial_{x} \lambda_{-n}^{+} +  \sum_{j=0}^{n} \lambda_{-j}^{+} \lambda_{j-n}^{+} + \sum_{m=1}^{n+2} \sum_{j=-1}^{n-m+1}  \frac{(-i)^m}{m!} \partial_{\sigma_{\perp}}^{m}  \lambda_{-j}^{+} \partial_{x_{\perp}}^{m} \lambda_{j-n+m}^{+} \nonumber \\
& \qquad +  \sum_{\substack{ j,l \\ \beta_{j} + \beta_{l} = - n } }  \lambda_{\beta_j}^{+} \lambda_{\beta_l}^{+} + \sum_{m=1}^{n+2} \sum_{\substack{ j,l \\ \beta_{j} + \beta_{l} = - n }} \frac{(-i)^m}{m!} \partial_{\sigma_{\perp}}^{m}  \lambda_{\beta_{j}}^{+} \partial_{x_{\perp}}^{m} \lambda_{\beta_{l-m}}^{+} 
 \Bigg). \label{Eqn.RecurInt}
\end{align}

Similarly, for the recursive terms of fractional order, we collect the terms of order $\beta_{n} = \alpha - n$ from \eqref{Eqn.OpMain3} to get
\begin{align} 
& \lambda_{\beta_{n+1}}^{-} \lambda_{1}^{+} + \lambda_{1}^{-} \lambda_{\beta_{n+1}}^{+}  
+ \sum_{j=0}^{n} \lambda_{\beta_{j}}^{-} \lambda_{j-n}^{+} + \sum_{j=0}^{n} \lambda_{j-n}^{-} \lambda_{\beta_{j}}^{+} \nonumber \\
& \qquad + \sum_{m=1}^{n+2} \sum_{j=0}^{n-m+1}  \frac{(-i)^m}{m!} \partial_{\sigma_{\perp}}^{m}  \lambda_{\beta_j}^{-} \partial_{x_{\perp}}^{m} \lambda_{j-n+m}^{+} = \partial_{x} \lambda_{\beta_{n}}^{+} \nonumber 
\end{align}
and again using \eqref{Eqn.010} we arrive at
\begin{align} 
\lambda_{\beta_{n+1}}^{\pm}  
&= \mp \frac{1}{2 \lambda_{1}^{+} } \Bigg( \partial_{x} \lambda_{\beta_{n}}^{+} + 2 \sum_{j=0}^{n} \lambda_{\beta_{j}}^{+} \lambda_{j-n}^{+}  \nonumber \\  & \qquad \qquad + \sum_{m=1}^{n+2} \sum_{j=0}^{n-m+1}  \frac{(-i)^m}{m!} \partial_{\sigma_{\perp}}^{m}  \lambda_{\beta_j}^{+} \partial_{x_{\perp}}^{m} \lambda_{j-n+m}^{+}   \Bigg). \label{Eqn.RecurFrac}
\end{align}

\section{Sweeping method for the Helmholtz equation}
\label{Section.Sweeping}

In this section we define a sweeping scheme to solve wave propagation problems for a solution $w$ with a time-harmonic dependence of the form $e^{- i \omega t}$. Wave propagation problems governed by \eqref{Eqn.LaplacianLocal}, for a time-harmonic source and time-harmonic boundary conditions with fixed frequency $\omega \in \mathbb{R}$, reduce to solving the Helmholtz equation
\begin{align} 
\HH w = \Delta w + \left(\omega^2 c^{-2} + i \omega a - a_{\alpha} (-i \omega)^{\alpha} \right) w = f.\label{Eqn.Helm}
\end{align}
The time dependence $e^{- i \omega t}$ is taken into account by replacing the time-derivatives in the wave operator \eqref{Eqn.LaplacianLocal}
by $-i\omega$. Equivalently, the terms $i \omega$ appearing in the symbolic version of pseudodifferential factorization \eqref{Eqn.Decomp01} are replaced by $- i \omega$, including the approximate symbols \eqref{Eqn.PadeSqrt}-\eqref{Eqn.Lambdabeta0} defined below. As a result, one obtains a pseudodifferential factorization for the Helmholtz equation \eqref{Eqn.Helm} where the frequency $\omega$ can be regarded as a fixed parameter. 

The above equation is considered in the half-space $\Omega = \{ x > 0 \} \subset \mathbb{R}^d$, augmented by the Sommerfeld radiation condition at infinity to guarantee a unique outgoing solution. The source $f$ is a compactly supported function. The compactly supported Dirichlet profile, denoted by $w_{\rm D}$, is imposed at the hyperplane $\{ x=0 \}$ . As before, we assume that the $x$-axis represents the dominant direction of wave propagation. We propose a double-sweep scheme based on the following decoupled system of marching equations,
\begin{align} 
\partial_{x} v - \Lambda_{M}^{-} v &= f, \label{Eqn.Sweep0} \\
\partial_{x} u - \Lambda_{M}^{+} u &= v, \label{Eqn.Sweep}
\end{align} 
where
\begin{align} 
\Lambda_{M}^{\pm} = \Op{ \lambda_{1,M}^{\pm} + \lambda_{0,M}^{\pm} + \lambda_{\beta_{1},M}^{\pm} }
\label{Eqn.TruncatedOp}
\end{align} 
and $\lambda_{n,M}^{\pm}$ is a Pad\'e approximation of order $M$ for the symbol $\lambda_{n}^{\pm}$ as defined below. Since $f$ is compactly supported, there is $0<L<\infty$ large enough for the half-space $\{ x \geq  L \}$ to be outside of the support of $f$. We impose the condition that $v=0$ in $\{ x \geq L \}$ as a boundary condition for $v$. Since $v=0$ in $\{ x \geq L \}$, then $\partial_{x} u = \Lambda_{M}^{+} u$ in $\{ x \geq L \}$ which means that $u$ is outgoing as required by the Sommerfeld radiation condition. We also impose $u = w_{\rm D}$ at $\{x = 0\}$ as the physical boundary condition.

Since we are including the highest degree symbols of the DtN maps $\Lambda^{\pm}$, we expect the accuracy of this sweeping method to increase as the frequency $\omega \to \infty$. In fact, the first term neglected from the expansion \eqref{Eqn.SymbolExpansion01} has a negative degree. This property is explored in greater detail in Section \ref{Section.ErrorAnalysis}. In order to obtain practical approximations for the operators associated with the symbols $\lambda_{1}^{\pm}$, $ \lambda_{0}^{\pm}$, and $\lambda_{\beta_1}^{\pm}$, we employ a high order, or so-called wide-angle, Pad\'e approximation for $\lambda_{1}^{\pm}$ in \eqref{Eqn.Sym1} as well as its appearance in the definitions of $ \lambda_{0}^{\pm}$, and $\lambda_{\beta_1}^{\pm}$ in \eqref{Eqn.Sym0} and \eqref{Eqn.SymB0}, respectively. This process renders a wide-angle approximation for these symbols and corresponding operators. Throughout, we will approximate functions of the type $(1 + z)^{\gamma}$ in the vicinity of $z = 0$, where $z = - |\sigma_{\perp}|^2 / (\omega^2 c^{-2})$. The underlying assumptions for these approximations to be accurate is that $|z| \leq \delta$ for some $\delta < 1$ \cite{BakerBook1996}, or equivalently that 
\begin{align} \label{eqn.delta}
\max (c) \, |\sigma_{\perp}| / \omega \leq \delta < 1.
\end{align}
Now we proceed to define the high order Pad\'e approximations for the symbols $\lambda_{1}^{\pm}$, $ \lambda_{0}^{\pm}$, and $\lambda_{\beta_1}^{\pm}$.

\textbf{Symbol of degree $\mathbf{1}$}: For the leading symbol $\lambda_{1}^{\pm}$ given by \eqref{Eqn.Sym1}, we use a high order Pad\'e approximation, in partial fraction form, for the function $(1+z)^{1/2}$ to obtain,
\begin{align} 
\lambda_{1,M}^{\pm} = \mp i \omega c^{-1} \left( a_{0}^{(1/2)} - \sum_{m=1}^{M} \frac{ \omega^2 c^{-2} a_{m}^{(1/2)} }{ \omega^2 c^{-2} b_{m}^{(1/2)} + |\sigma_{\perp}|^2 } \right) \label{Eqn.PadeSqrt}
\end{align}
where the coefficients $a_{m}^{(1/2)}$ and $b_{m}^{(1/2)}$ are described in more detail in the Appendix \ref{Section.Pade}. 

\textbf{Symbol of degree $\mathbf{0}$}: For the symbol $\lambda_{0}^{\pm}$ given by \eqref{Eqn.Sym0}, we apply a Pad\'e approximation of the function $(1+z)^{-1/2}$ for the second term to arrive at
\begin{align} \label{Eqn.Lambda0comp}
\lambda_{0,M}^{\pm} &= \partial_{x} \ln \left( \left( i \lambda_{1,M}^{+} \right)^{\mp 1/2} \right) \mp \frac{a c}{2} \left( a_{0}^{(-1/2)} - \sum_{m=1}^{M} \frac{ \omega^2 c^{-2} a_{m}^{(-1/2)} }{ \omega^2 c^{-2} b_{m}^{(-1/2)} + |\sigma_{\perp}|^2 } \right) \\
& \qquad \pm \frac{i c}{2 \omega} \left( \nabla_{x_{\perp}} \ln c \right) \cdot (i \sigma_{\perp}) \left( a_{0}^{(-3/2)} - \sum_{m=1}^{M} \frac{ \omega^2 c^{-2} a_{m}^{(-3/2)} }{ \omega^2 c^{-2} b_{m}^{(-3/2)} + |\sigma_{\perp}|^2 } \right) \nonumber
\end{align}
where the Pad\'e coefficients $a_{m}^{(-1/2)}$ and $b_{m}^{(-1/2)}$ and $a_{m}^{(-3/2)}$ and $b_{m}^{(-3/2)}$, are described in the Appendix \ref{Section.Pade} and $(i \sigma_{\perp})$ is the symbol of the tangential derivatives. The first term in the right-hand side of \eqref{Eqn.Lambda0comp} is conveniently written to be integrated analytically. For this purpose, it is necessary to make use of the Pad\'e approximation of the functions  $(1+z)^{1/4}$ and $(1+z)^{-1/4}$ to obtain
\begin{align} 
& \sqrt{ \frac{i \lambda_{1,M}^{+}}{\omega c^{-1}} } =   a_{0}^{(1/4)} - \sum_{m=1}^{M} \frac{ \omega^2 c^{-2} a_{m}^{(1/4)} }{ \omega^2 c^{-2} b_{m}^{(1/4)} + |\sigma_{\perp}|^2 } 
\label{Eqn.Lambda0+1-4} \\
& \sqrt{ \frac{\omega c^{-1}}{i \lambda_{1,M}^{+} } } =  a_{0}^{(-1/4)} - \sum_{m=1}^{M} \frac{ \omega^2 c^{-2} a_{m}^{(-1/4)}}{ \omega^2 c^{-2} b_{m}^{(-1/4)} + |\sigma_{\perp}|^2 } 
\label{Eqn.Lambda0-1-4}
\end{align}
where the Pad\'e coefficients $a_{m}^{(1/4)}$ and $b_{m}^{(1/4)}$, $a_{m}^{-(1/4)}$ and $b_{m}^{-(1/4)}$, are described in the Appendix \ref{Section.Pade}.

\textbf{Symbol of degree $\boldsymbol{\beta_{1}}$}: For the fractional attenuation symbol $\lambda^{\pm}_{\beta_1}$ given by \eqref{Eqn.SymB0}, we employ a Pad\'e approximation for the the function $(1+z)^{-1/2}$ to obtain
\begin{align} 
\lambda_{\beta_1,M}^{\pm} = \mp \frac{a_{\alpha} (i \omega)^\alpha}{2 i \omega c^{-1}} \left( a_{0}^{(-1/2)} - \sum_{m=1}^{M} \frac{ \omega^2 c^{-2} a_{m}^{(-1/2)} }{ \omega^2 c^{-2} b_{m}^{(-1/2)} + |\sigma_{\perp}|^2 } \right) \label{Eqn.Lambdabeta0}
\end{align}
where the Pad\'e coefficients $a_{m}^{(-1/2)}$ and $b_{m}^{(-1/2)}$ are described in the Appendix \ref{Section.Pade}.

For these three symbols, we have the following behavior for the difference between the exact symbols and the Pad\'e approximations,
\begin{align} 
&| \lambda_{1,M}^{\pm} - \lambda_{1}^{\pm} | = \mathcal{O} \left( \omega \delta^{4M+2} \right),  \label{Eqn.401} \\
&| \lambda_{0,M}^{\pm} - \lambda_{0}^{\pm} | = \mathcal{O} \left( \delta^{4M+2} \right), \label{Eqn.402} \\
&| \lambda_{\beta_1,M}^{\pm} - \lambda_{\beta_1}^{\pm} | = \mathcal{O} \left( \omega^{\alpha - 1} \delta^{4M+2} \right).  \label{Eqn.403}
\end{align}
provided that \eqref{eqn.delta} is satisfied.

We conclude this section by noticing that in order to solve \eqref{Eqn.Sweep0}-\eqref{Eqn.Sweep}, where $\Lambda_{M}^{\pm}$ is given by \eqref{Eqn.TruncatedOp}, and the symbols $\lambda_{1,M}^{\pm}$, $\lambda_{0,M}^{\pm}$, and $\lambda_{\beta_1,M}^{\pm}$ are expressed in Pad\'e fraction forms as in \eqref{Eqn.PadeSqrt}, \eqref{Eqn.Lambda0comp} and \eqref{Eqn.Lambdabeta0}, respectively, then it is required to obtain the operator associated with symbols of the form $a/b$ where $a=a(x_{\perp},\omega)$ is independent of time $t$ and of the spatial frequency $\sigma_{\perp}$, and $b=b(x_{\perp},\omega,\sigma_{\perp})$ is independent of time $t$. Hence, following details shown in Section \ref{Section.ReviewPDO}, we have that $\Op{a/b} \sim \Op{a} \Op{1/b}$. The operator of the numerator $a$ is straightforward to compute. The operator of the inverse of the denominator $b$ is more involved. Following Section \ref{Section.ReviewPDO}, we have that
\begin{align} 
\Op{1/b} = \Op{b}^{-1} + \mathcal{R} \nonumber
\end{align} 
where $\Op{b}$ is an operator of order $+2$, $\Op{b}^{-1}$ is an operator of order $-2$, and $\mathcal{R}$ is an operator of order $-3$, provided that the inverse of $\Op{b}$ exists. In all the cases included in this section, the denominators of the Pad\'e fractions are invertible because the ($\theta$-rotated) Pad\'e coefficients possess an imaginary part as explained in Appendix \ref{Section.Pade}. In the rest of this paper, we neglect the remainders $\mathcal{R}$.

\section{Modeling error analysis}
\label{Section.ErrorAnalysis}

Here we analyze the extent to which the solution to the proposed pseudodifferential sweeping method satisfies the Helmholtz equation. This analysis of residuals is then the first step to estimate the difference between the solution $w$ to the Helmholtz equation \eqref{Eqn.Helm} and solution $u$ to the proposed pseudodifferential sweeping method \eqref{Eqn.Sweep0}-\eqref{Eqn.Sweep}. 
By construction $w$ and $u$ satisfy the same Dirichlet boundary values in the hyperplane $\{x=0\}$. So it will be useful to consider the following problem in the half-space. Let $\Omega = \{(x,x_{\perp}) \in \mathbb{R}^d : x>0, x_{\perp} \in \mathbb{R}^{d-1} \}$, and consider the following boundary value problem for the Helmholtz operator \eqref{Eqn.Helm},
\begin{equation} \label{Helmholtz_eq}
\HH v = f \quad \text{in $\Omega$}, \quad v=0 \quad \text{on} \quad \partial \Omega = \{x=0\},
\end{equation}
satisfying the Sommerfeld radiation condition at infinity, and where $c(x,x_{\perp})\geq c_0>0$, $a(x,x_{\perp})\geq a_0>0$ and $a_\alpha(x,x_{\perp}) \geq 0$ are smooth functions. Here it is assumed that $c = \text{const.}$, $a=\text{const.}$, $a_{\alpha} = 0$, and $f=0$ outside of a bounded subdomain of $\Omega$. The existence, uniqueness, and regularity of solutions for this problem can be established (in strong and weak formulations) through the method of images, ie., by extending $c$, $a$, and $a_{\alpha}$ symmetrically about $\{x=0\}$, and $f$ anti-symmetrically about $\{x=0\}$. Thus, by posing the problem in the whole-space $\mathbb{R}^d$, the classical results for well-posedness can be invoked \cite{McLean2000,Lio-Mag-Book-1972}. Adding to the well-posedness, we now review the behavior of the solution with respect to the frequency $\omega$.

\begin{lemma} \label{Thm.Lemma1}
Under the assumption on the domain $\Omega$, the wave speed $c$, the damping coefficient $a$ and the fractional attenuation coefficient $a_{\alpha}$ stated above, there is a constant $C>0$ independent of $\omega$, such that for all $f\in L^2(\Omega)$, 
\begin{align} \label{Eqn.C1}
\sum_{k=0}^{2} \omega^{1-k} \|v\|_{H^{k}(\Omega)} \leq C\|f\|_{L^2(\Omega)}
\end{align}
where $v$ is a solution to \eqref{Helmholtz_eq}. The same holds for the (formal) adjoint $\mathcal{H}^*$.
\end{lemma}

\begin{proof}
From the weak formulation of the Helmholtz equation, we obtain
\begin{align*} \label{Eqn.C2}
& - \|\nabla v\|^2_{L^2(\Omega)} +\omega^2\|c^{-1} v\|^2_{L^2(\Omega)} + i\omega \|\sqrt{a} v\|^2_{L^2(\Omega)} - (-i\omega)^\alpha\|\sqrt{a_\alpha} v\|^2_{L^2(\Omega)} = \int_\Omega f v \, dx
\end{align*}
where the integration over the unbounded domain $\Omega$ is finite due to the exponential decay of the solution for non-vanishing damping ($a \geq a_{0} > 0$). See details in \cite{ColtonKressBook2013,McLean2000}.  From the imaginary part of the previous equality, we deduce 
\begin{align*} 
\omega \|\sqrt{a} v\|^2_{L^2(\Omega)} + \sin (\pi \alpha/2) \, \omega^{\alpha} \|\sqrt{a_{\alpha}} v\|^2_{L^2(\Omega)} \leq \|f\|_{L^2(\Omega)} \|v\|_{L^2(\Omega)}
\end{align*}
and where $\alpha \in (0,1)$ and $\text{Im} (- (-i \omega)^{\alpha}) = \omega^{\alpha} \text{Im}( - e^{-i\pi \alpha / 2}) = \omega^{\alpha} \sin (\pi \alpha /2) > 0 $, therefore, 
\begin{align*} \label{Eqn.C3}
\omega \|v\|_{L^2(\Omega)} \leq C \|f\|_{L^2(\Omega)}.
\end{align*}

On the other hand, the real part of the equality implies $\|\nabla v\|^2_{L^2(\Omega)} \leq C (\omega^2 + \omega^{\alpha}\cos(\pi \alpha / 2)) \|v\|^2_{L^2(\Omega)} + \|f\|_{L^2(\Omega)} \|v\|_{L^2(\Omega)}$, which combined with the previous inequality renders
\begin{align*} \label{Eqn.C4}
\|v\|_{H^1(\Omega)}\leq C\|f\|_{L^2(\Omega)}
\end{align*}
for a constant $C$ independent of $\omega$, for all $\omega \geq 1$ so that $\omega^2 \geq \omega^{\alpha}$. Now, due to elliptic regularity, for a source $f \in L^2(\Omega)$, the solution $v \in H^{2}(\Omega)$ satisfying \eqref{Helmholtz_eq} strongly. Therefore, we can take the $L^2(\Omega)$-norm of \eqref{Helmholtz_eq} and combine with previous inequalities to find that $\| \Delta v \|_{L^2(\Omega)} \leq C \omega \| f \|_{L^2(\Omega)}$ where again $C$ is independent of $\omega$ for all sufficiently large $\omega$. 
A standard argument based on finite quotients and the bound on $\Delta v$ (see \cite{Evans2008,McLean2000} for details) leads to the improved regularity estimate
\begin{align*} \label{Eqn.C5}
\|v\|_{H^2(\Omega)}\leq C \omega \|f\|_{L^2(\Omega)},
\end{align*}
for another constant $C$ independent of $\omega$. The proof concludes by adding the three estimates above.
\end{proof}

\begin{lemma} \label{Thm.Lemma2}
Under the assumptions of Lemma \ref{Thm.Lemma1}, there is a constant $C>0$ independent of $\omega$ such that, for all $f\in H^{-2}(\Omega)$, 
\begin{align} \label{Eqn.C6}
\|v\|_{L^2(\Omega)}\leq C \omega \|f\|_{H^{-2}(\Omega)}
\end{align}
where $v$ is the solution (by transposition) of \eqref{Helmholtz_eq}.
\end{lemma}

\begin{proof}
By definition, $v$ satisfies
\begin{align*} \label{Eqn.C7}
\langle v,\mathcal{H}^*\phi\rangle = \langle f, \phi\rangle
\end{align*}
for all $\phi\in H^2(\Omega)\cap H^1_0(\Omega)$. For an arbitrary $\psi\in L^2(\Omega)$ we take $\phi$ solution to $\mathcal{H}^{*}\phi = \psi$ with boundary condition $\phi = 0$ as in Lemma \ref{Thm.Lemma1}. Then, 
\begin{align*} \label{Eqn.C8}
|\langle v,\psi\rangle| &\leq \|f\|_{H^{-2}(\Omega)} \|\phi\|_{H^2(\Omega)} \leq C \omega \|f\|_{H^{-2}(\Omega)} \|\psi\|_{L^2(\Omega)},
\end{align*}
for the same constant $C>0$ from Lemma \ref{Thm.Lemma1}, which renders the desired inequality.
\end{proof}

With the above results we can quantify the behavior for the difference between $w$ and $u$ as the frequency $\omega$ increases. 
Like in the previous section, this analysis is intended for wave fields that are highly oscillatory in time (with dependence $e^{-i\omega t}$) and also in space provided that the waves travel in a relatively narrow neighborhood of the $x$-axis. Such solutions to the wave equation are mathematically described by being supported in phase space in the set $\{ \varrho \, \omega^2 \leq \sigma^2 + |\sigma_{\perp}|^2 \}$ intersected with the wedge $\{ c \, |\sigma_{\perp}| \leq \delta \, \omega \}$ for some constants $\varrho > 0$ and $0 < \delta < 1$. These conditions are illustrated in Figure \ref{Fig:Diagram2} and stated in precise terms in Assumption \ref{supp_fourier}. As shown below, this parameter $\delta$ quantifies the degree to which the waves propagate primarily along the $x$-axis, allowing the principal symbol $\lambda_{1}^{\pm}$ to remain elliptic, defining the type of wave field solutions for which we can expect accuracy using the sweeping method \eqref{Eqn.Sweep0}-\eqref{Eqn.Sweep}, and how much more accuracy can be gained by increasing the order $M$ of the Pad\'e approximation.

\begin{figure}[htbp]
\centering
\includegraphics[width=0.65 \textwidth, trim = 5 10 5 10, clip]{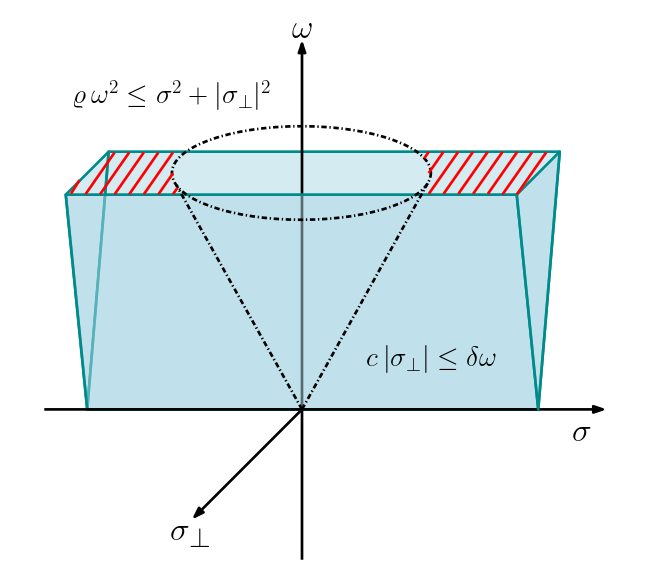}
\caption{Support in Fourier space (shaded in red) of the solution $u$ in order for the sweeping method \eqref{Eqn.Sweep0}-\eqref{Eqn.Sweep} to be accurate. These conditions are stated in Assumption \ref{supp_fourier} and employed in Theorem \ref{Thm.theorem1}. In physical terms, these conditions ensure that the waves, induced by the source $f$, oscillate in space at least as fast as they oscillate in time, and simultaneously that the $x$-axis is the dominant direction of propagation.}
\label{Fig:Diagram2}
\end{figure}

\begin{assumption} \label{supp_fourier} 
Assume that there exist constants $0 < \delta < 1$ and $0 < \varrho$ independent of $\omega$ and $M$ such that the Fourier transform $\hat{u}$ of the solutions $u$ to the proposed pseudodifferential sweeping method \eqref{Eqn.Sweep0}-\eqref{Eqn.Sweep} satisfies $ \text{supp} (\hat{u}) \subset  \{ (\sigma , \sigma_{\perp}) \in \mathbb{R}^d \,:\, \varrho \, \omega^2 \leq \sigma^2 + |\sigma_{\perp}|^2 \,  \text{and} \,  \max(c) | \sigma_{\perp}| \leq \delta \, \omega  \} $.
\end{assumption}

Recall the Helmholtz operator $\HH$ from \eqref{Eqn.Helm} and the boundary value problem \eqref{Helmholtz_eq}. Let $\HH^{-1} : H^{-2}(\Omega) \to H^{0}(\Omega)$ denote the operator that provides the solution to the Helmholtz equation in $\Omega$ given a prescribed source $f$ and a homogeneous Dirichlet boundary condition at the boundary of $\Omega$. See \cite[Ch. 2 \S 6]{Lio-Mag-Book-1972} for details on well-posedness of elliptic equations in Sobolev scales including solutions by transposition. The $\HH^{-1}$ operator is well defined for all $\omega > 0$ due to the Sommerfeld radiation condition, non-vanishing damping ($a \geq a_{0} > 0$) and non-negative fractional attenuation ($a_{\alpha} \geq 0$) coefficients.

Now, define the residual
\begin{align} \label{Eqn.Residual}
r = f - \HH u = \HH \left( w - u \right).
\end{align}
For non-vanishing damping, the error can be estimated from the residual as follows
\begin{align} 
\|w - u \|_{H^{0}(\Omega)} = \| \HH^{-1} r \|_{H^{0}(\Omega)} \leq C \omega \| r \|_{H^{-2}(\Omega)} \label{Eqn.Error}
\end{align}
with a constant $C=C(\Omega,c,a, a_{\alpha}, \alpha) > 0$ independent of $\omega$ as shown by Lemma \ref{Thm.Lemma1} and Lemma \ref{Thm.Lemma2}. Now, for our proposed pseudodifferential method based on \eqref{Eqn.Sweep0}-\eqref{Eqn.Sweep} and \eqref{Eqn.PadeSqrt}-\eqref{Eqn.Lambdabeta0}, the exact DtN maps $\Lambda^{\pm}$ are approximated by $\Lambda_{M}^{\pm}$ satisfying
\begin{align} \label{Eqn.Res01}
\Lambda^{\pm} = \Lambda_{M}^{\pm} + \RR_{M}^{\pm}
\end{align}
where $M$ denotes the number of Pad\'e terms, the remainder operators $\RR_{M}^{\pm} \in \Op{\SS^1}$ whose symbols $\rho_{M}^{\pm}$ satisfy
\begin{align} \label{Eqn.Res02}
\rho_{M}^{\pm} = \rho_{\rm pade}^{\pm} + \rho_{\rm trunc}^{\pm}
\end{align}
where $\rho_{\rm pade}^{\pm} \in \SS^1$ is due to the Pad\'e approximation of the symbols $\lambda_{1}^{\pm}$, $\lambda_{0}^{\pm}$ and $\lambda_{\beta_1}^{\pm}$ described in the previous section, and $\rho_{\rm trunc}^{\pm} \in \SS^{-1}$ is due to the truncation \eqref{Eqn.TruncatedOp} of the classical expansion \eqref{Eqn.SymbolExpansion01}. As a consequence of \eqref{Eqn.401}-\eqref{Eqn.403} and the first neglected term in \eqref{Eqn.SymbolExpansion01}, we have
\begin{align} 
\rho_{\rm Pade}^{\pm} &\sim \omega \delta^{4M+2} + \delta^{4M+2} + \omega^{\alpha-1} \delta^{4M+2} \lesssim  \omega \, \delta^{4M},
\label{Eqn.Res05} \\
\rho_{\rm trunc}^{\pm} &\sim \omega^{-1}.\label{Eqn.Res06}
\end{align}

Now, the residual satisfies
\begin{align} \label{Eqn.Res08}
r &= f - \HH u  = f -  ( \partial_{x} - \Lambda^{-})( \partial_{x} - \Lambda^{+}) u  \nonumber \\
&= f - ( \partial_{x} - \Lambda_{M}^{-} - \RR_{M}^{-})( \partial_{x} - \Lambda_{M}^{+} - \RR_{M}^{+}) u  \nonumber \\
&= f - ( \partial_{x} - \Lambda_{M}^{-})( \partial_{x} - \Lambda_{M}^{+}) u + ( \partial_{x} - \Lambda_{M}^{-}) \RR_{M}^{+} u \nonumber \\ 
& \qquad + \RR_{M}^{-} ( \partial_{x} - \Lambda_{M}^{+}) u - \RR_{M}^{-}\RR_{M}^{+} u \nonumber \\
&= ( \partial_{x} - \Lambda_{M}^{-}) \RR_{M}^{+} u + \RR_{M}^{-} ( \partial_{x} - \Lambda_{M}^{+}) u - \RR_{M}^{-}\RR_{M}^{+} u,
\end{align}
and we make use of the equivalence of the Sobolev norms defined in Fourier space (see \eqref{Eqn.SobNormPsiOp} in Section \ref{Section.ReviewPDO}) to compute the $H^{-2}$ norm of $r$ in terms of the behavior of the symbols for $\Lambda_{M}^{\pm}$ and  $\RR_{M}^{\pm}$, as follows,
\begin{align} \label{Eqn.Res09}
 \| r \|^2_{H^{-2}(\Omega)} &\leq C \int \left( \frac{\left( |\sigma| + \omega \right)^2 \left( \omega \delta^{4M} + \omega^{-1} \right)^2 }{\left(1 + \sigma^2 + |\sigma_{\perp}|^2\right)^2} + \frac{ \left( \omega \delta^{4M} + \omega^{-1}  \right)^4 }{ \left(1 + \sigma^2 + |\sigma_{\perp}|^2\right)^2 } 
 \right) |\hat{u}|^2 d (\sigma, \sigma_{\perp}) \nonumber \\
& \leq C \left( \delta^{4M} + \omega^{-2} \right)^2 \| u \|^2_{H^{0}(\Omega)}
\end{align}
where we have used the following inequality
\begin{align} 
\frac{\left( |\sigma| + \omega \right)^2 }{\left(1 + \sigma^2 + |\sigma_{\perp}|^2\right)^2} &\leq  \frac{ 2 (\sigma^2 + \omega^2) }{ \left(1 + \sigma^2 + |\sigma_{\perp}|^2\right)^2 } \nonumber \\ 
&\leq  \frac{2}{\left(1 + \sigma^2 + |\sigma_{\perp}|^2\right)}  + \frac{2 \omega^2}{\left(1 + \sigma^2 + |\sigma_{\perp}|^2\right)^2} \nonumber \\
& \leq \frac{2}{\varrho \omega^2} + \frac{2 \omega^2}{\varrho^2 \omega^4} \leq \frac{C}{\omega^2} 
\nonumber
\end{align}
which is valid under the conditions of Assumption \ref{supp_fourier}. The constant $C>0$ in \eqref{Eqn.Res09} is independent of all sufficiently large frequencies $\omega$, $\delta < 1$ and Pad\'e order $M \geq 1$. The above arguments have shown the following result.

\begin{theorem} \label{Thm.theorem1}
Under the conditions of Lemma \ref{Thm.Lemma1} and Assumption \ref{supp_fourier}, the following bound
\begin{align} \label{Eqn.Res003}
\frac{ \| u - w \|_{H^{0}(\Omega)} }{ \| u \|_{ H^{0}(\Omega)} } \leq C \left( \omega \delta^{4M} + \frac{1}{\omega} \right)
\end{align}
holds for solutions $w$ to the Helmholtz equation \eqref{Eqn.Helm} and corresponding solutions $u$ to the proposed pseudodifferential sweeping method \eqref{Eqn.Sweep0}-\eqref{Eqn.Sweep}, for all sufficiently large frequencies $\omega$, any number of Pad\'e terms $M \geq 1$, where the constant $C>0$ is independent of $\omega$, $\delta$ and $M$.
\end{theorem}

As a result, we observe that the error $\| u - w \|_{H^{0}(\Omega)}$ can be controlled if the number $M$ of Pad\'e terms increases as the frequency $\omega$ of oscillations also increases. More specifically, since $0 < \delta < 1$, then $\delta = e^{- \beta}$ for some $\beta > 0$. If $M=M(\omega) \sim \ln \omega /(2 \beta)$, then $\delta^{4M} \sim 1/\omega^2$ and thus the overall error decreases as $1/\omega$. These estimates are tested numerically in the next section.

\section{Numerical experiments}
\label{Section.Numerical}

In this section we provide the results from some numerical experiments to illustrate the implementation of the proposed pseudodifferential sweeping method. For computational purposes, we are forced to consider a bounded domain and apply absorbing boundary conditions or layers to approximate the effect of the Sommerfeld radiation condition.

We consider a square domain $\Omega \subset \mathbb{R}^2$ of side $L=1$ centered at the origin. At the center of the domain, there is a circular inclusion of radius $0.1$ with wavespeed $c=2$ while the rest of the domain has a background wavespeed $c_{o} = 1$. However, a smooth transition is accomplished by the following definition,
\begin{align} \label{Eqn.WS}
c(x,y) = 1 + H( 0.1 - \sqrt{ x^2 + y^2 })
\end{align}
where $H(s) = 1/(1+e^{-800s})$ is a smooth version of the Heaviside function. The (unattenuated) wave number is $k(x,y) = \omega / c(x,y)$.

\begin{figure}[htbp]
\centering
\captionsetup{font=small}
\subfloat[Wavespeed $c$]{
\includegraphics[width=0.35 \textwidth, trim = 0 0 0 0, clip]{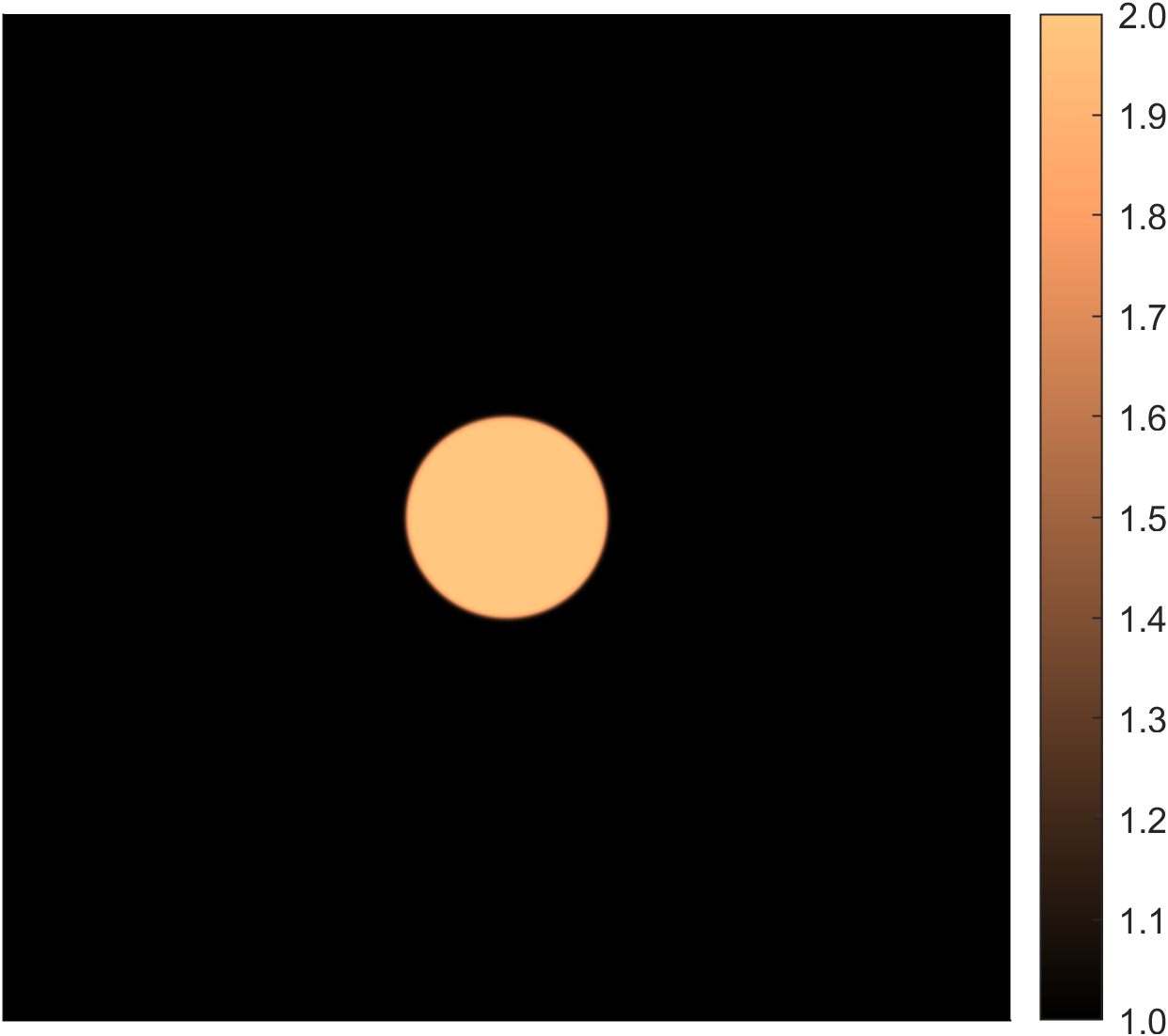}
\label{subfig:wavespeed}
}
\subfloat[Absorbing sponge $\beta$]{
\includegraphics[width=0.35 \textwidth, trim = 0 0 0 0, clip]{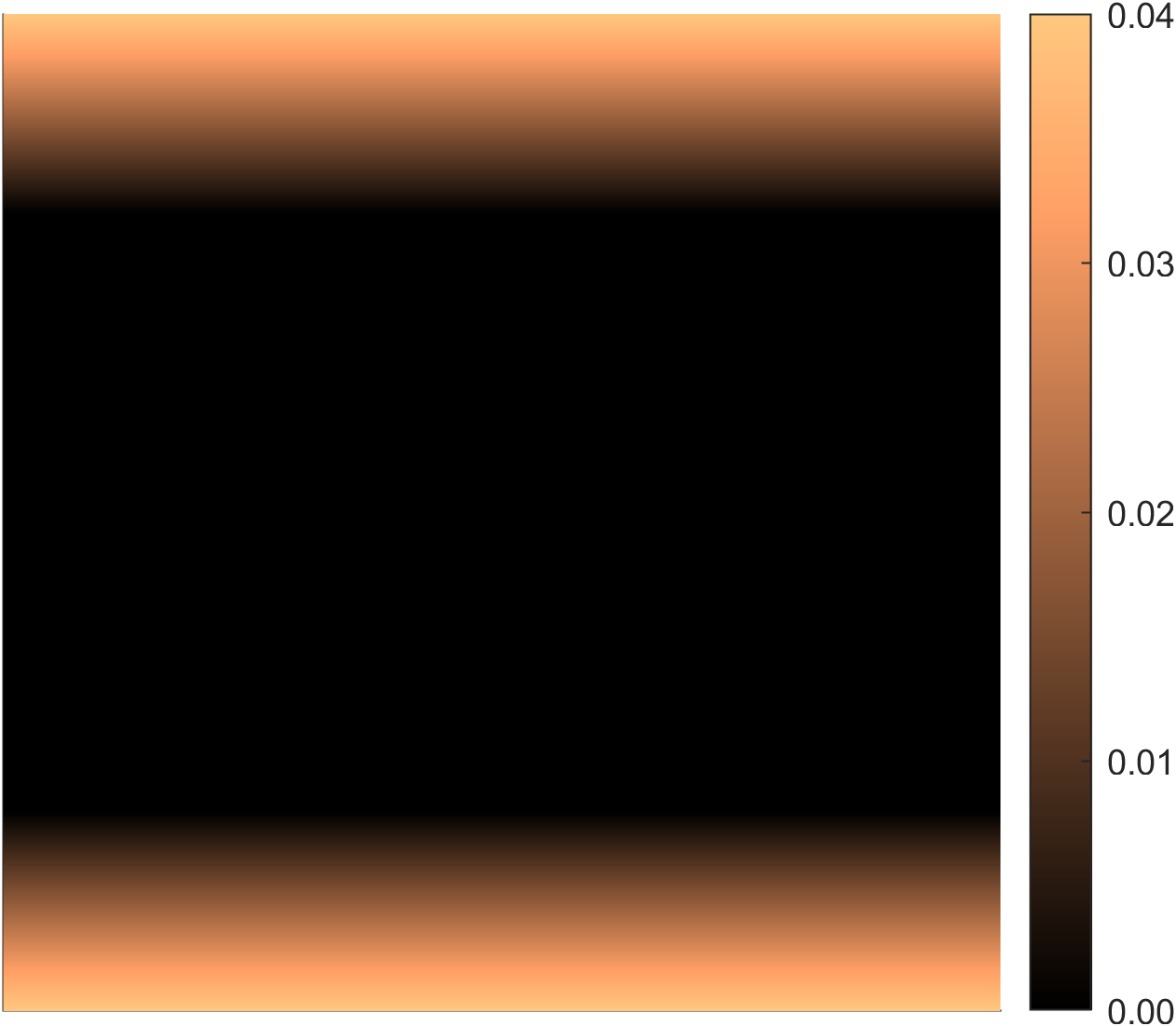}
\label{subfig:absorbing}
} 
\caption{(a) Wavespeed profile with an inclusion defined by \eqref{Eqn.WS} and (b) absorbing sponge defined by \eqref{Eqn.Sponge} to mitigate effects of the top and bottom boundaries.}
\label{Fig.Medium_PML}
\end{figure}

On the left boundary, a Dirichlet boundary profile
\begin{align} 
u(y) = e^{- i \omega /c_{o} L/2 } e^{- 200 y^2} \label{Eqn.BCleft}
\end{align}
is imposed. On the right boundary, an outgoing condition is imposed by setting $v=0$. Also, in order to mitigate the influence of the top and bottom boundaries, a sponge boundary layer is introduced where the wavenumber $k$ is replaced by $k (1 + i \beta(x,y))/\sqrt{ 1 + \beta(x,y)^2} $ with
\begin{align} \label{Eqn.Sponge}
\beta(x,y) = 
\begin{cases}
    0.2 (  |y| - 0.3  ), & \text{if $|y| > 0.3$}\\
    0, & \text{otherwise}.
\end{cases}
\end{align}
As a result, an exponential decay of the solution is observed inside the sponge layer (in the vicinity of the top and bottom boundaries). The actual boundary conditions at the top and bottom boundaries are first order absorbing condition of the type $\partial_{n} u =  i k u$ where $\partial_{n}$ denotes the derivative in the outward normal direction, and $k$ is the modified version of the wavenumber. Similar approaches have been used in \cite{Stolk2017}. The wave speed and sponge layer are illustrated in Figure \ref{Fig.Medium_PML}. In all the experiments, the damping coefficient is set to $a=0.01$, the fractional attenuation coefficient $a_{\alpha} = 10$ and the fractional exponent $\alpha = 0.5$.

The sweeping method is implemented in two steps. First, a \emph{one-way solution} $u_{\rm one}$ is computed by solving 
\begin{align}
\partial_{x} u_{\rm one} - \Op{ \lambda_{1,M}^{+} + \lambda_{0,M}^{+} + \lambda_{\beta_1,M}^{+}  } u_{\rm one} = 0. \label{Eqn.701one}
\end{align}
and the prescribed boundary condition \eqref{Eqn.BCleft} on the left boundary. Then, the \emph{two-way solution} $u_{\rm two}$ is defined as $u_{\rm two} = u_{\rm one} + u$ where $u$ solves
\begin{align}
\partial_{x} v - \Op{ \lambda_{1,M}^{-} + \lambda_{0,M}^{-} + \lambda_{\beta_1,M}^{-}  } v &= - \HH u_{\rm one}, \label{Eqn.703two} \\
\partial_{x} u - \Op{ \lambda_{1,M}^{+} + \lambda_{0,M}^{+} + \lambda_{\beta_1,M}^{+}  } u &= v, \label{Eqn.705two}
\end{align}
for vanishing Dirichlet condition for $v$ on the right boundary, and vanishing Dirichlet condition for $u$ on the left boundary. Notice that both solutions $u_{\rm one}$ and $u_{\rm two}$ satisfy the prescribed boundary condition \eqref{Eqn.BCleft} on the left boundary and the outgoing boundary condition at the right boundary by virtue of \eqref{Eqn.701one} and \eqref{Eqn.705two} with $v=0$ on the right boundary. However, as opposed to $u_{\rm one}$, we expect $u_{\rm two}$ to account for reflection effects from the wavespeed inclusion in the domain $\Omega$ due to the incorporation of the source $- \HH u_{\rm one}$ in \eqref{Eqn.703two}.

For the numerical results presented here, \eqref{Eqn.701one}-\eqref{Eqn.705two} were discretized using Heun's method for the stepping in the $x$-direction using $36$ points per wavelength. All of the second-order $y$-derivatives corresponding to the symbol $|\sigma_{\perp}|^2$ appearing in \eqref{Eqn.PadeSqrt}-\eqref{Eqn.Lambdabeta0}
were discretized using a second-order centered finite difference scheme using $12$ points per wavelength. The pseudocode is shown in Algorithm \ref{Algo.1}. The
numerical approximations for the relative residuals $\| r \|_{H^{-2}(\Omega)} / \| u \|_{H^{0}(\Omega)}$ are displayed in Table \ref{tab:table1}. Several runs for 
doubling the frequency $\omega$ and linearly increasing the number of Pad\'e terms $M$ were performed. This configuration tests the estimates \eqref{Eqn.Res003} obtained in Section \ref{Section.ErrorAnalysis} for the error and residual to decrease as $1/\omega$ as the frequency $\omega$ increases and the number of Pad\'e terms grows logarithmically as $M \sim \log \omega$. For each row in the table and for the diagonal, the observed orders of decay were obtained as the slope of a straight line fitted through the numerical residuals versus frequency in the log-log space. These numerical results conform with the expected behavior derived in \eqref{Eqn.Res003}.

\begin{algorithm}[ht]
\caption{\label{Algo.1} Sweeping algorithm for a general one-way equation of the form $\partial_{x} u - \Op{\lambda} u = F$ using Heun's stepping scheme. The number of steps in the $x$-direction is $N_x$ and the step-size is $\Delta x = L/N_x$. Here $\lambda = \lambda_{1,M}^{\pm} + \lambda_{0,M}^{\pm}+ \lambda_{\beta_1,M}^{\pm}$ where these symbols are defined in Section \ref{Section.Sweeping} and depend on the sweeping variable $x$ through the dependence of the wave speed $c$, damping coefficient $a$ and fractional attenuation coefficients $a_{\alpha}$.} 
\begin{algorithmic} 
\STATE Initialize by applying a boundary condition: $u_{0} = w_{\rm D}$
\FOR {$i=1,2,\ldots, N_{x}$}
    \STATE Predictor:
    \STATE $\qquad \qquad p = \Op{ \lambda(x_{i-1}) } u_{i-1} + F_{i-1}$
    \STATE $\qquad \qquad u_{\rm aux} = u_{i-1} + \Delta x \, p$
    \STATE Corrector:
    \STATE $\qquad \qquad q = \Op{ \lambda(x_{i}) } u_{\rm aux} + F_{i}$
    \STATE $\qquad \qquad u_{i} = u_{i-1} + \Delta x \, (p + q)/2$
\ENDFOR
\end{algorithmic} 
\end{algorithm}

\begin{figure}[htbp]
\centering
\captionsetup{font=small}
\subfloat[Real part of one-way solution]{
\includegraphics[width=0.4 \textwidth, trim = 0 0 0 0, clip]{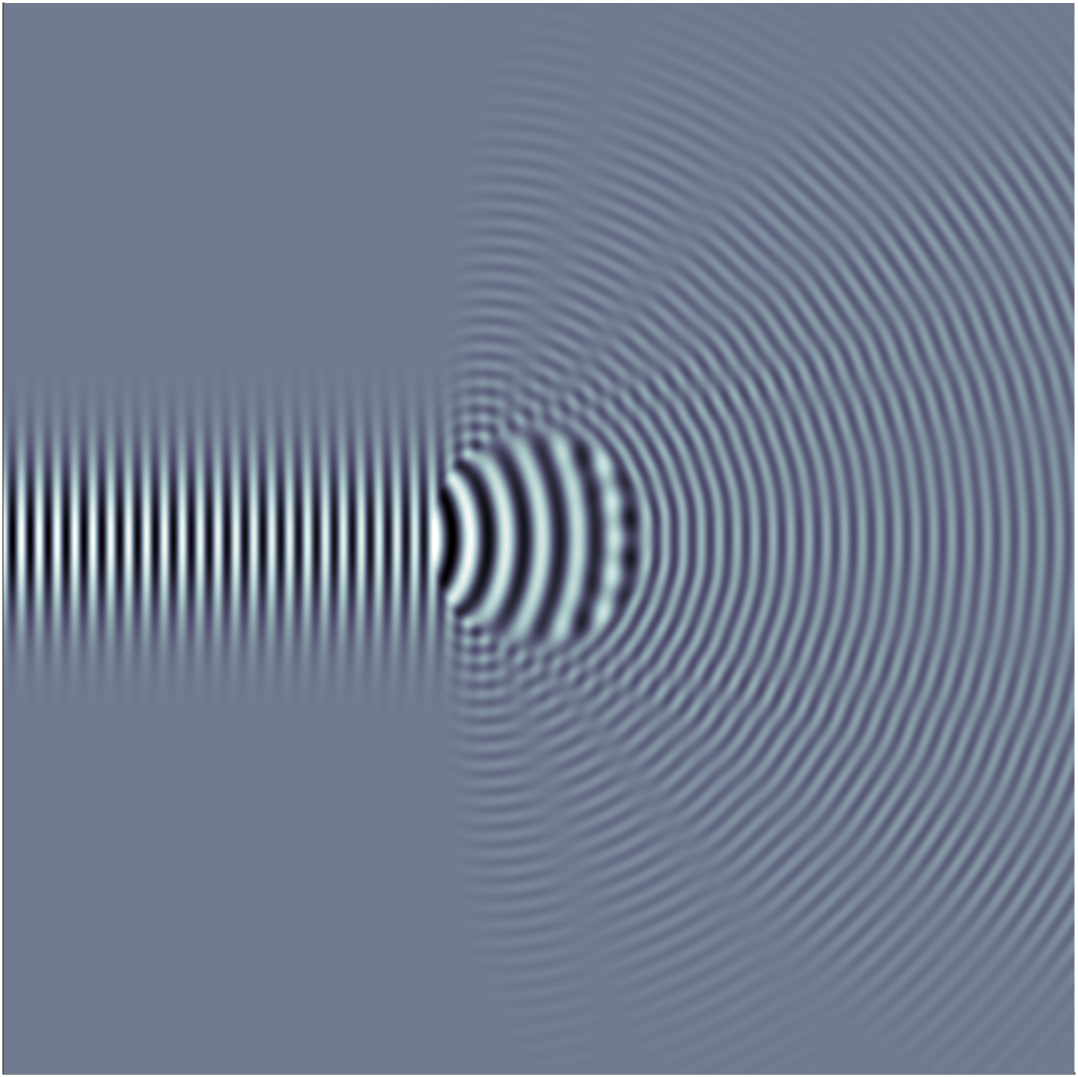}
\label{subfig:RealPart1}
}
\subfloat[Real part of two-way solution]{
\includegraphics[width=0.4 \textwidth, trim = 0 0 0 0, clip]{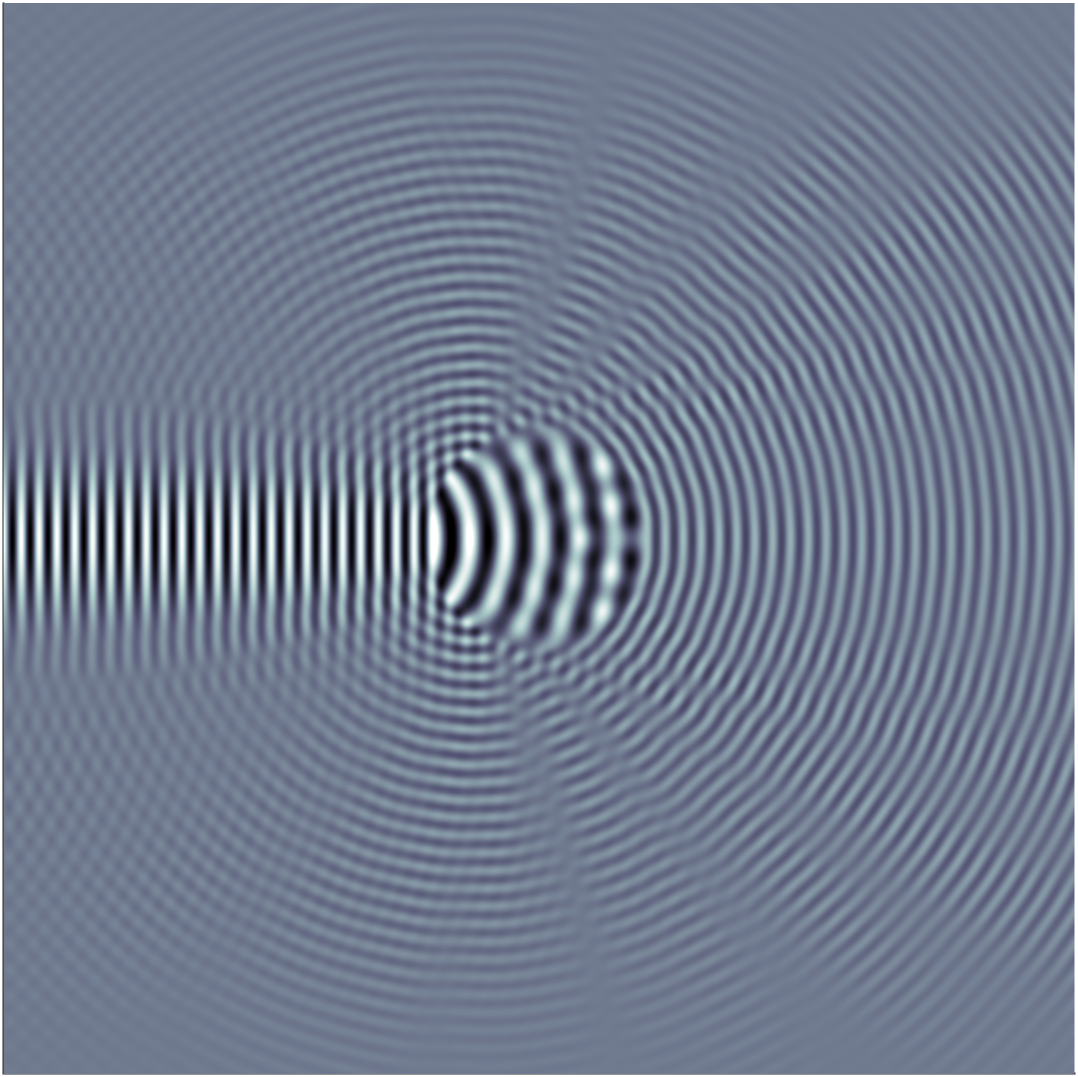}
\label{subfig:RealPart2}
} \\
\subfloat[Amplitude of one-way solution]{
\includegraphics[width=0.4 \textwidth, trim = 0 0 0 0, clip]{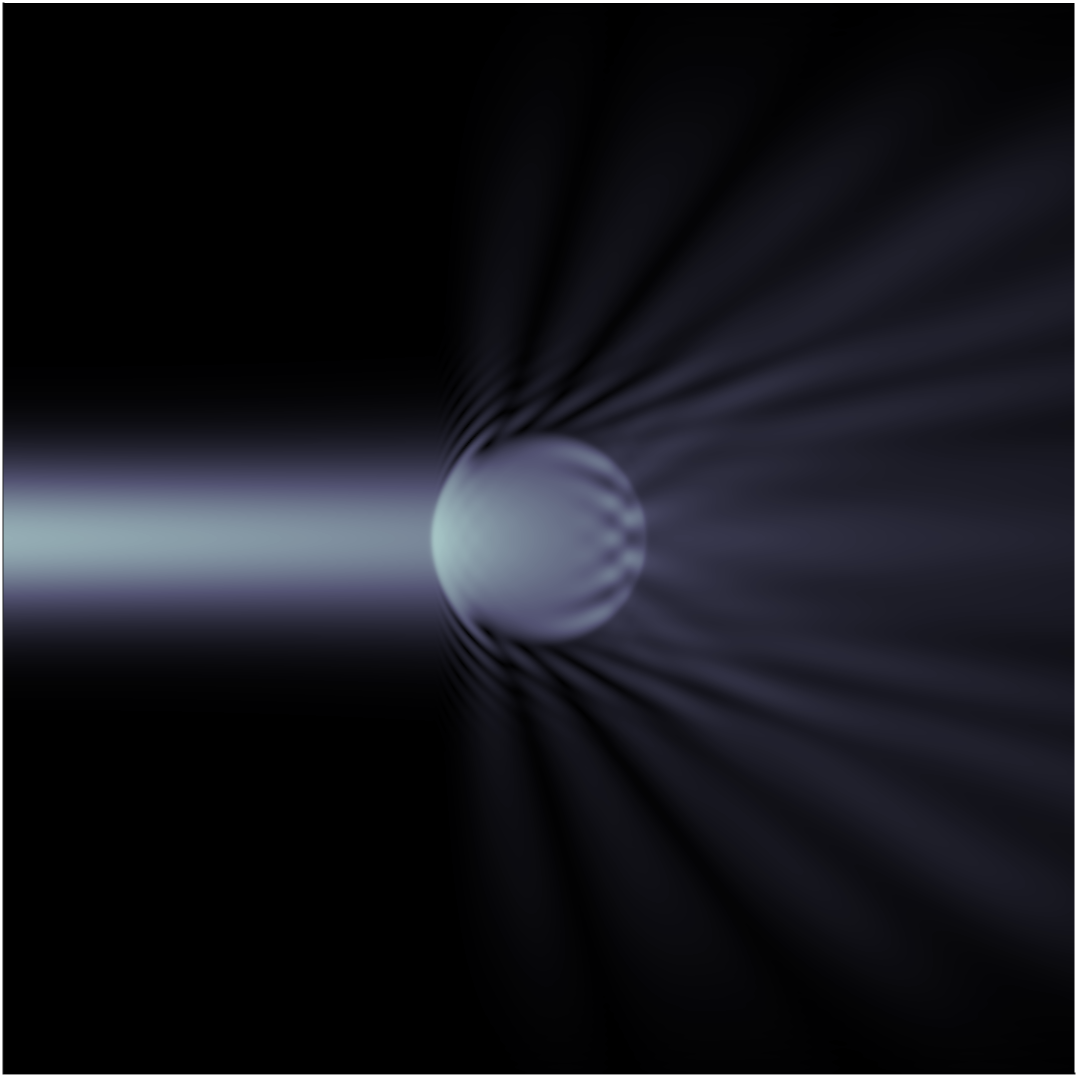}
\label{subfig:Amplitude1}
}
\subfloat[Amplitude of two-way solution]{
\includegraphics[width=0.4 \textwidth, trim = 0 0 0 0, clip]{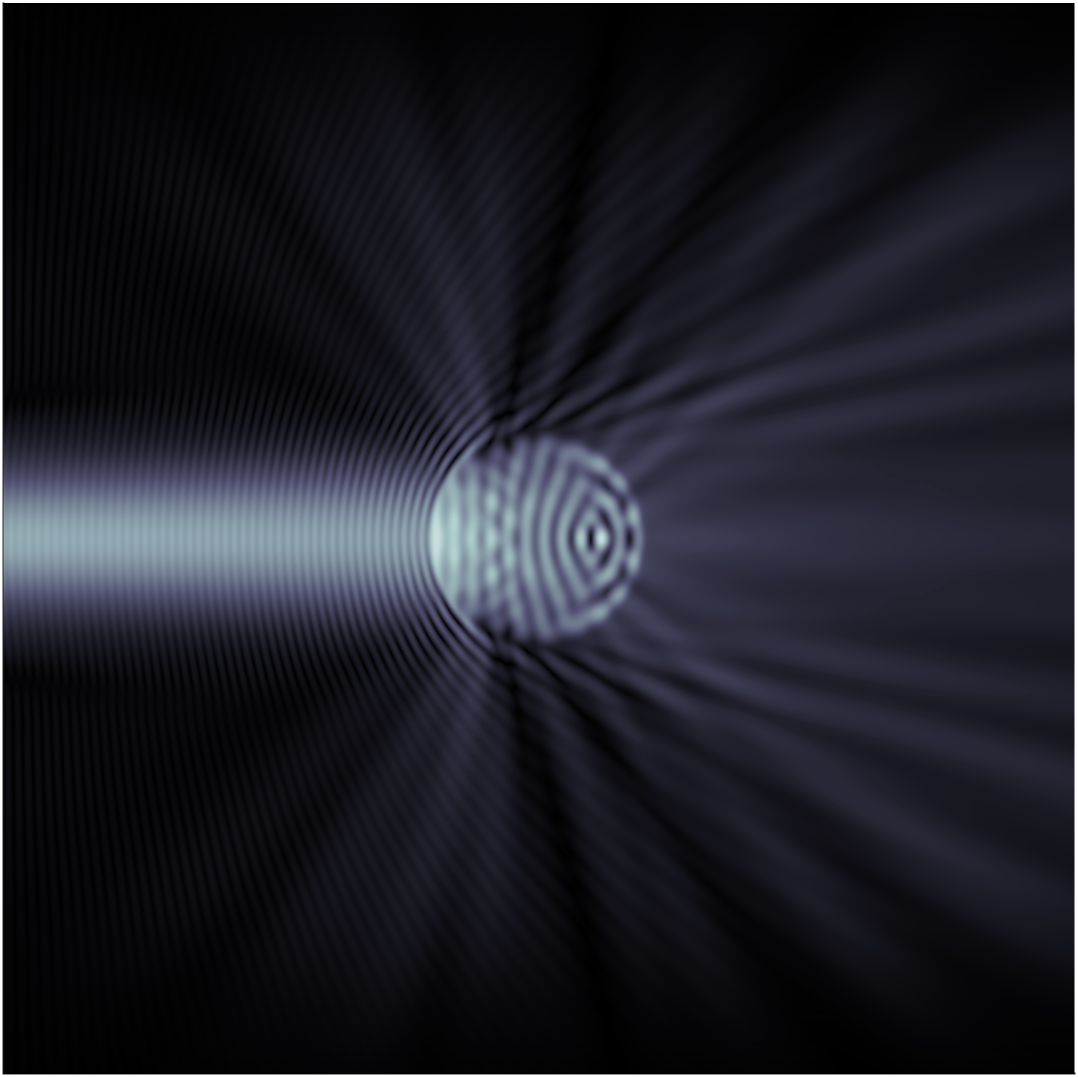}
\label{subfig:Amplitude2}
} \\
\subfloat[FFT of one-way solution]{
\includegraphics[width=0.4 \textwidth, trim = 0 0 0 0, clip]{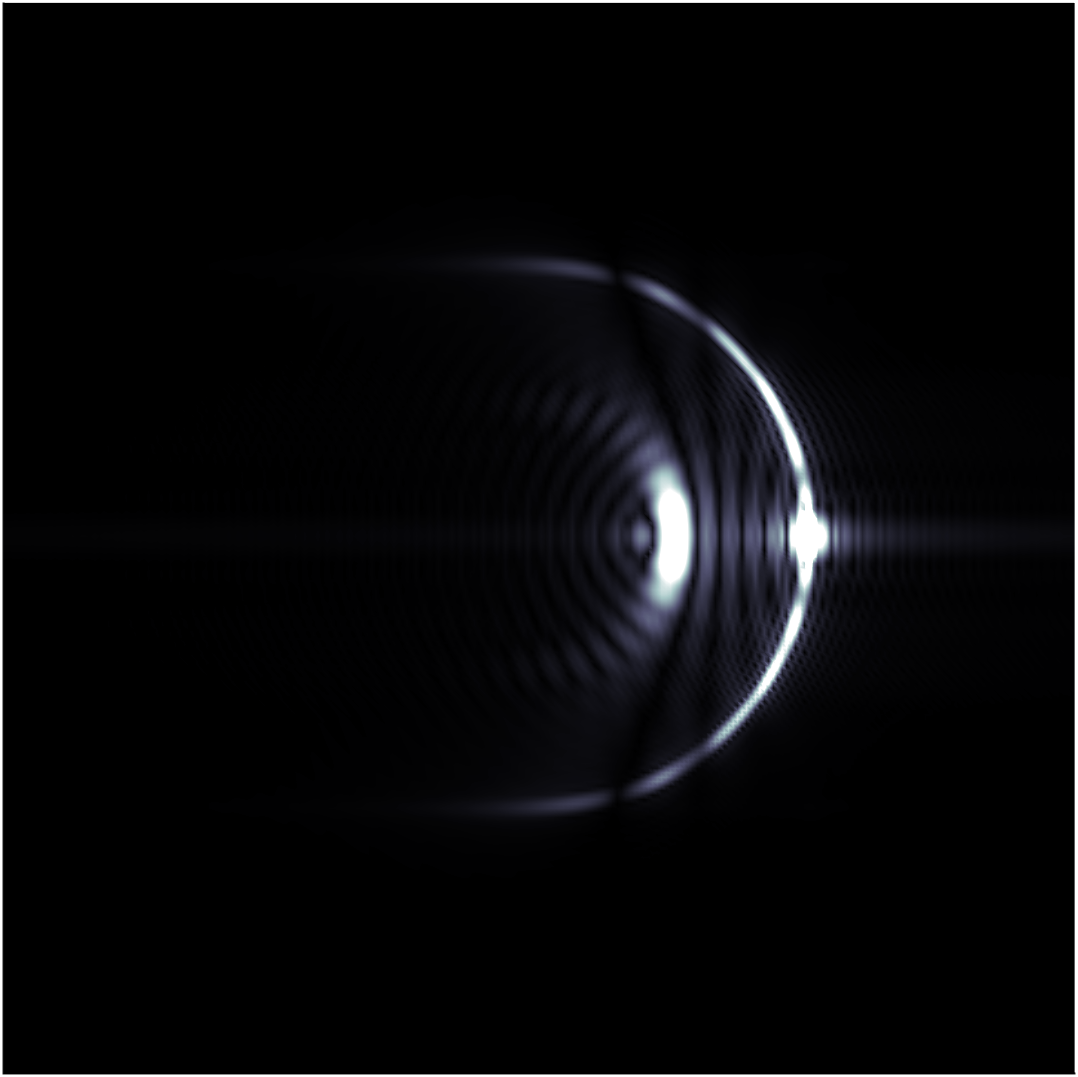}
\label{subfig:FFT1}
}
\subfloat[FFT of two-way solution]{
\includegraphics[width=0.4 \textwidth, trim = 0 0 0 0, clip]{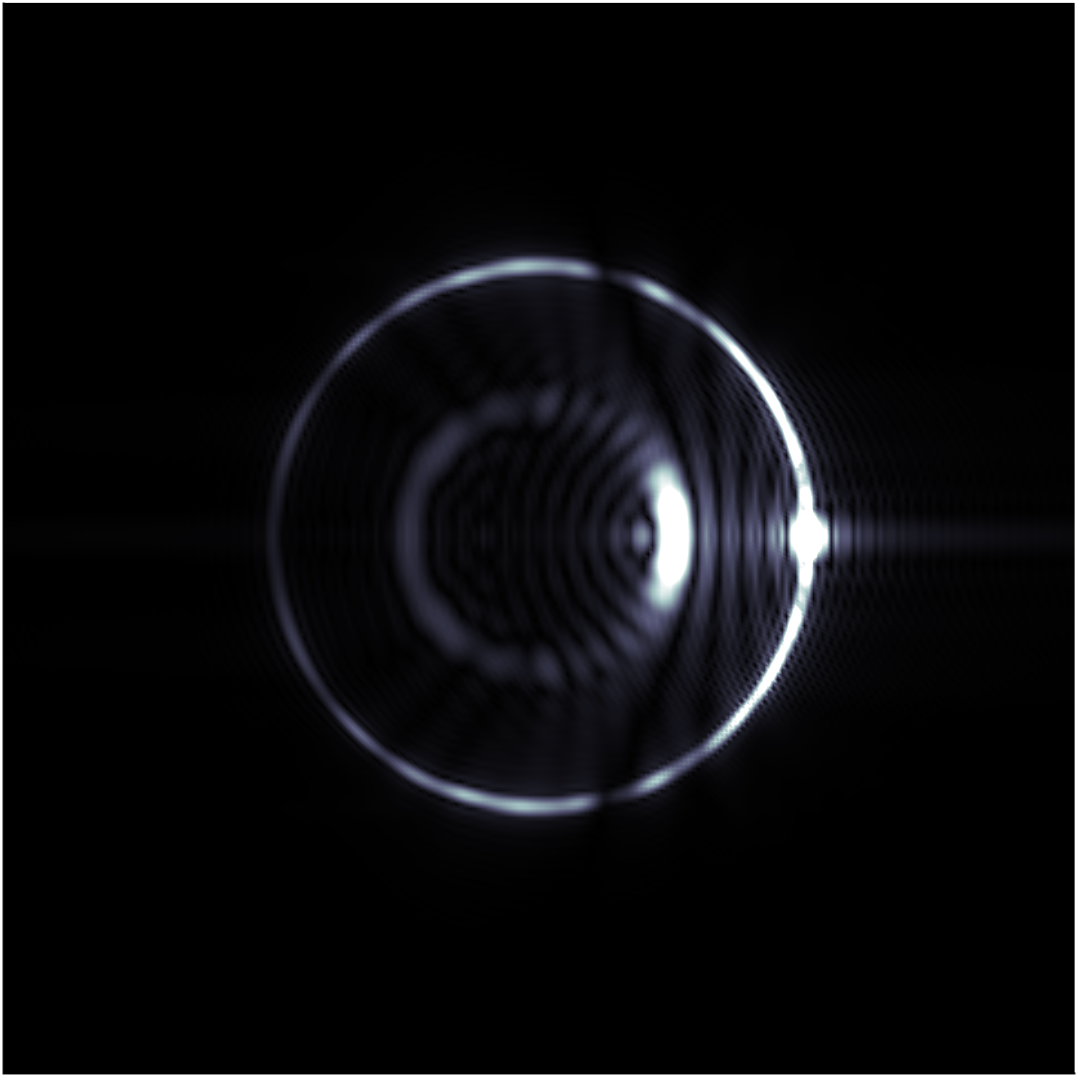}
\label{subfig:FFT2}
}
\caption{Comparison between the one-way and two-way numerical solutions. These solutions were computed using Pad\'e approximations of the pseudodifferential symbols with $4$ terms. The frequency is $\omega = 120 \pi$ which fits $60$ wavelengths across the domain. The two-way solution captures the reflections induced by the inclusion. These reflections are visible on plot of the real part and amplitude of the solution, as well as in the Fourier plots.}
\label{Fig.Solution_OMEGA=60}
\end{figure}

\begin{table}[h]
\centering
\caption {\label{tab:table1} Numerical approximations for the relative residuals $\| r \|_{H^{-2}(\Omega)} / \| u \|_{H^{0}(\Omega)}$ for doubling the frequency $\omega$ and linearly increasing the number of Pad\'e terms $M$. This configuration tests the estimates \eqref{Eqn.Res003} obtained in Section \ref{Section.ErrorAnalysis} for the error and residual to decrease as $1/\omega$ as the frequency $\omega$ increases and the number of Pad\'e terms grows logarithmically $\sim \log \omega$. The observed orders of decay were obtained as the slope of a straight line fitted through the numerical residuals versus frequency in the log-log space.} 
\begin{tabular}{lrrrrr}
\hline
\hline
Pad\'e terms & \multicolumn{4}{c}{Relative residual (\%)} & Observed order \\
\hline
 &  $\omega = 20 \pi$  & $\omega = 40 \pi$  & $\omega = 80 \pi$ & $\omega = 160 \pi$ &  \\  
\cline{2-5} 

$3$ & $2.73$ & $1.85$ & $1.53$ & $1.33$ & $-0.34$ \\
 
$4$ & $2.46$ & $1.30$ & $0.63$ & $0.33$ & $-0.97$ \\
 
$5$ & $2.54$ & $1.35$ & $0.65$ & $0.29$ & $-1.04$ \\

$6$ & $2.54$ & $1.35$ & $0.64$ & $0.29$ & $-1.04$ \\

\hline
& \multicolumn{4}{r}{Observed order along the diagonal} & $-1.07$ \\

\hline
\hline
\end{tabular}
\end{table}

\section{Conclusion and limitations}
\label{Section.Conclusion}

The Helmholtz equation is notoriously difficult to solve in part because the domains of dependence and influence of its solutions are global, i.e., the solution value at any specific point can affect and be affected by the solution at any other point. This phenomenon is due to the ability of waves to propagate over long distances. Hence, an effective solver (or preconditioner) must either account for such global behavior \cite{Ernst2012, Gander2019} or incorporate a-priori assumptions about the direction of wave propagation. The pseudodifferential method presented here is designed to do the latter effectively. The pseudodifferential calculus allows us to incorporate the physical effects of variable media properties. These effects include: 
\begin{enumerate}[leftmargin=2em,rightmargin=0em]
    \item Propagation and refraction modeled by the principal symbol $\lambda_{1}^{\pm}$ defined in \eqref{Eqn.Sym1}.
    \item Amplitude modulation modeled by the first and third terms of the symbol $\lambda_{0}^{\pm}$  defined in \eqref{Eqn.Sym0}.
    \item Damping modeled by the second term of the symbol $\lambda_{0}^{\pm}$ defined in \eqref{Eqn.Sym0}.
    \item Fractional attenuation modeled by the symbol $\lambda_{\beta_1}^{\pm}$ defined in \eqref{Eqn.SymB0}.
    \item Reflection modeled by solving both the forward and the backward propagation equations \eqref{Eqn.Sweep0}-\eqref{Eqn.Sweep}.
    \item Wide angle of wave propagation thanks to high order Pad\'e approximations applied to the square root symbol appearing not only in the principal symbol \eqref{Eqn.Sym1} but also in \eqref{Eqn.Sym0} and \eqref{Eqn.SymB0}.
\end{enumerate}
Moreover, since the neglected pseudodifferential terms decay as the frequency $\omega$ increases, the proposed sweeping method is well-suited for biomedical applications based on high-frequency ultrasonics. The authors are in the process of developing numerical implementations of the proposed pseudodifferential methods to obtain ultrasound tomographic imaging and run high-intensity focused ultrasound simulations. As soon as meaningful results are obtained from these efforts, they will be reported in forthcoming publications.  

In this work, we have included a proof-of-concept numerical implementation using Heun's stepping method and second-order finite difference discretization for tangential derivatives. The purpose of the numerical implementation was to test the theoretical error estimates. However, we did not make a detailed analysis of numerical error and computational cost associated with the proposed discretization. The use of higher-order stepping schemes (Runge-Kutta or exponential integrators) and efficient approximations of the tangential derivatives (spectral, finite element, or compact finite difference methods), and their computational advantages should be the subject of future studies. The interplay between the stepping scheme, the Pad\'e approximations of the pseudodifferential symbols, and the step sizes on the $x$-axis and tangential directions has an impact on the numerical stability of the method. Hence, CFL-type conditions should be established for each scheme and grid refinement.

In general, the proposed pseudodifferential method is limited by its underlying assumptions, namely, that waves propagate along a dominant direction and that the media properties are smooth. Hence, we expect this method to be less accurate when waves encounter steep variations in material properties that induce refraction away from the dominant direction of propagation. Discontinuous jumps in media properties could be handled accurately by properly incorporating transmission conditions into the forward and backward sweeps. However, such approach remains to be developed and tested. 

\appendix
\section{Pad\'e approximations} \label{Section.Pade}

The real-valued Pad\'e approximant for the function $(1+z)^{\gamma}$, as a partial fraction expansion, has the following form,
\begin{align} \label{Eqn.PadeAppendix}
P_{M}(z) =  a_{0} + \sum_{m=1}^{M} \frac{a_{m} }{ z - b_{m} }
\end{align}
where the real-valued coefficients $a_{m}$ and $b_{m}$ are computed in order to match the value of the function $(1+z)^{\gamma}$ and its first $2M$ derivatives at $z=0$ \cite{BakerBook1996}. Unfortunately, these Pad\'e approximants suffer from inaccuracies and instabilities when $z \leq -1$ due to the branch cut along the negative real line starting at $z=-1$ \cite{Lu1998a,Milinazzo1997,Modave2020}. In order to avoid this problem (for the specific case of $\gamma=1/2$) Milinazzo et al. \cite{Milinazzo1997} considered the Pad\'e approximation of the function $e^{i \theta \gamma} \left( 1 + \zeta \right)^{\gamma}$ where $\zeta = (1+z) e^{-i \theta} - 1$, which has a rotated branch cut defined by the angle $\theta$. As a result, the approximation remains stable and continuous for $z \in \mathbb{R}$.  This $\theta$-rotated Pad\'e approximant of order $M$, in partial fraction form, is given by
\begin{align} \label{Eqn.PadeRotatedAppendix}
P_{M,\theta}(z) &= e^{i \theta \gamma} \left( a_{0} + \sum_{m=1}^{M} \frac{a_{m} }{ (1+z) e^{- i \theta} - 1 - b_{m} } \right) = \tilde{a}_{0} + \sum_{m=1}^{M} \frac{\tilde{a}_{m} }{ z - \tilde{b}_{m} }
\end{align}
where the complex-valued Pad\'e coefficients are $\tilde{a}_{0} = a_{0} e^{i \theta \gamma }$, $\tilde{a}_{m} = a_{m} e^{i \theta (1+\gamma)}$ and $\tilde{b}_{m} = (1+b_{m})e^{i \theta} - 1$ for $m=1,2,...,M$ \cite{BakerBook1996,Milinazzo1997}.
In Section \ref{Section.Sweeping}, there is need to apply the Pad\'e method for the cases $\gamma = \pm 1/2$ to approximate the symbols $(\lambda_{1}^{+})^{\pm 1}$, and the cases $\gamma = \pm 1/4$ for the symbols $(\lambda_{1}^{+})^{\pm 1/2}$. For these cases, the real-valued Pad\'e coefficients, displayed on Tables \ref{Tbl.01}-\ref{Tbl.04}, were computed using Sanguigno's \textsc{MATLAB} \emph{pade} function to obtain a rational approximation \cite{SanguignoPade2023}, followed by R2020b \textsc{MATLAB}'s \emph{residue} function to re-cast it as a partial fraction expansion. The well-known error analysis for Pad\'e approximants (see for instance \cite{BakerBook1996}) shows that
\begin{align} \label{Eqn.PadeError}
(1+z)^{\gamma} = P_{M,\theta}(z)  + \mathcal{O} \left( |z|^{2M + 1} \right)
\end{align}
as $z \to 0$, which justifies the error terms specified in Section \ref{Section.Sweeping} by replacing $z = - c^{2} \delta^2$.

\begin{table}[ht]
\scriptsize
\centering
\caption{Real-valued Pad\'e coefficients to approximate the function $(1+z)^{1/2}$ in partial fraction form \eqref{Eqn.PadeAppendix} for orders $M=1,2,3,4$.}
\begin{tabular}{l | lllll | lllll}
\hline
M & $a_{0}$  & $a_{1}$ & $a_{2}$ & $a_{3}$ & $a_{4}$  & $b_{1}$ & $b_{2}$ & $b_{3}$ & $b_{4}$   \\
\hline
1 & 2.8889 & -7.1358 &  & & & -3.7778 & & &   \\
2  & 4.7738 & -34.5138 & -0.2786 & & & -9.6264 & -1.4778 & &  \\
3 & 6.7228 & -98.1129 & -1.0233 & -0.0723 & & -18.7042 & -2.4499 & -1.2139 & \\
4 & 8.6939 & -213.677 & -2.4055 & -0.2410 & -0.0303 & -31.0166 & -3.8025 & -1.6590 & -1.1242 \\
\hline
\end{tabular} \label{Tbl.01}
\end{table}

\begin{table}[ht]
\scriptsize
\centering
\caption{Real-valued Pad\'e coefficients to approximate the function $(1+z)^{-1/2}$ in partial fraction form \eqref{Eqn.PadeAppendix} for orders $M=1,2,3,4$.}
\begin{tabular}{l | lllll | lllll}
\hline
M & $a_{0}$  & $a_{1}$ & $a_{2}$ & $a_{3}$ & $a_{4}$  & $b_{1}$ & $b_{2}$ & $b_{3}$ & $b_{4}$   \\
\hline
1 & 0.3590 & 0.8218 &  & & & -1.2821 & & &   \\
2  & 0.2114 & 1.0955 & 0.4181 & & & -2.6945 & -1.0944 & &  \\
3 & 0.1493 & 1.4521 & 0.4474 & 0.2885 & & -4.9532 & -1.5846 & -1.0480 & \\
4 & 0.1152 & 1.8313 & 0.5184 & 0.2862 & 0.2219 & -8.0266 & -2.3254 & -1.3120 & -1.0292 \\
\hline
\end{tabular} \label{Tbl.02}
\end{table}

\begin{table}[ht]
\scriptsize
\centering
\caption{Real-valued Pad\'e coefficients to approximate the function $(1+z)^{1/4}$ in partial fraction form \eqref{Eqn.PadeAppendix} for orders $M=1,2,3,4$.}
\begin{tabular}{l | lllll | lllll}
\hline
M & $a_{0}$  & $a_{1}$ & $a_{2}$ & $a_{3}$ & $a_{4}$  & $b_{1}$ & $b_{2}$ & $b_{3}$ & $b_{4}$   \\
\hline
1 & 1.6239 & -1.5572 &  & & & -2.4957 & & &   \\
2  & 2.0906 & -5.8925 & -0.1638 & & & -6.0834 & -1.3433 & &  \\
3 & 2.4805 & -14.0723 & -0.4809 & -0.0572 & & -11.6531 & -2.1510 & -1.1613 & \\
4 & 2.8202 & -26.8939 & -0.9694 & -0.1539 & -0.0284 & -19.2030 & -3.2899 & -1.5524 & -1.0952 \\
\hline
\end{tabular} \label{Tbl.03}
\end{table}

\begin{table}[ht]
\scriptsize
\centering
\caption{Real-valued Pad\'e coefficients to approximate the function $(1+z)^{-1/4}$ in partial fraction form \eqref{Eqn.PadeAppendix} for orders $M=1,2,3,4$.}
\begin{tabular}{l | lllll | lllll}
\hline
M & $a_{0}$  & $a_{1}$ & $a_{2}$ & $a_{3}$ & $a_{4}$  & $b_{1}$ & $b_{2}$ & $b_{3}$ & $b_{4}$   \\
\hline
1 & 0.6213 & 0.5737 &  & & & -1.5149 & & &   \\
2  & 0.4794 & 1.1826 & 0.1945 & & & -3.3516 & -1.1593 & &  \\
3 & 0.4035 & 1.9393 & 0.3201 & 0.1095 & & -6.2432 & -1.7354 & -1.0796 & \\
4 & 0.3547 & 2.8260 & 0.4555 & 0.1670 & 0.0734 & -10.1697 & -2.5810 & -1.3816 & -1.0481 \\
\hline
\end{tabular} \label{Tbl.04}
\end{table}

\vspace{2em}

\section*{Acknowledgments}
The authors would like to thank the anonymous referees for their helpful suggestions that improved the quality of the manuscript considerably.

\vspace{2em}


\bibliographystyle{siamplain}
\bibliography{library}
\end{document}

%% file: ex_shared.tex
\usepackage{lipsum}
\usepackage{amsfonts}
\usepackage{graphicx}
\usepackage{epstopdf}
\usepackage{algorithmic}
\ifpdf
  \DeclareGraphicsExtensions{.eps,.pdf,.png,.jpg}
\else
  \DeclareGraphicsExtensions{.eps}
\fi

\newsiamremark{remark}{Remark}
\newsiamremark{hypothesis}{Hypothesis}
\crefname{hypothesis}{Hypothesis}{Hypotheses}
\newsiamthm{claim}{Claim}

\headers{Pseudodifferential Models for Ultrasound Waves}{S Acosta, J Chan, R Johnson and B Palacios}

\title{Pseudodifferential Models for Ultrasound Waves with Fractional Attenuation\thanks{Submitted to the editors DATE.
\funding{The work of B. Palacios was partially supported by Agencia Nacional de Investigaci\'on y Desarrollo (ANID), Grant FONDECYT Iniciaci\'on N$^\circ$11220772. J. Chan and R. Johnson gratefully acknowledge funding from the NSF through awards DMS-1943186 and DMS-2231482. S. Acosta would like to thank the support provided by Texas Children’s Hospital.}}}

\author{Sebastian Acosta\thanks{Department of Pediatrics, Baylor College of Medicine and Texas Children's Hospital, Houston, TX, USA 
  (\email{sacosta@bcm.edu}).}
\and Jesse Chan\thanks{Department of Computational Applied Mathematics and Operations Research, Rice University, Houston, TX, USA.}
\and Raven Johnson\footnotemark[3]
\and Benjamin Palacios\thanks{Department of Mathematics, Pontificia Universidad Catolica de Chile, Santiago, Chile.} }

\usepackage{amsopn}

\def\Op#1{ \text{Op}\left( #1 \right)}
\def\Sym#1{ \text{Sym}\left( #1 \right)}

\def\HH{\mathcal{H}}
\def\LL{\mathcal{L}}

\def\RR{\mathcal{R}}

\def\SS{\mathcal{S}}

\newtheorem{assumption}[theorem]{\textit{Assumption}}

\usepackage{amssymb}
\usepackage{amsbsy}
\usepackage{hyperref}
\usepackage{subfig}
\usepackage{enumitem}
